\documentclass[11pt]{article}

\usepackage{lineno,hyperref}
\usepackage{rotating}
\usepackage{makeidx}
\usepackage{float}
\usepackage{epsfig}
\usepackage{amsmath,amssymb}
\usepackage[latin1]{inputenc}
\usepackage{epic,eepic}
\usepackage{overpic}
\usepackage{color}
\usepackage{amsfonts}
\usepackage{graphicx}
\usepackage{mwe} % For dummy images
\usepackage{subcaption}
\usepackage{enumitem}

\usepackage{algorithm}
\usepackage{algorithmicx}
\usepackage{algpseudocode}

\newenvironment{demo}{\smallskip\noindent{{\it Proof.}}\hskip \labelsep}%
            {\hfill\penalty10000\raisebox{-.09em}{\large\bf\rm $\blacksquare$}\par\medskip}

\newtheorem{theorem}{Theorem}[section]
\newtheorem{lemma}[theorem]{Lemma}
\newtheorem{proposition}[theorem]{Proposition}

\begin{document}
\makeatletter{\renewcommand*{\@makefnmark}{}

%\small

%\begin{frontmatter}

\title{\bf{Non-linear WENO B-spline based approximation method}}

\footnotetext{Dr. Sergio Amat and Dr. Juan Ruiz have been supported by the national research project PID2019-108336GB-I00. Dr. Juan Ruiz has been also supported by the Fundaci\'on S\'eneca  grant  21728/EE/22  (Este  trabajo  es  resultado  de  la  estancia  (21728/EE/22)  financiada por  la  Fundaci\'on S\'eneca-Agencia de Ciencia y Tecnolog\'ia de la Regi\'on de Murcia con cargo al Programa Regional de Movilidad, Colaboraci\'on Internacional e Intercambio de Conocimiento ``Jim\'enez de la Espada'' (Plan de Actuaci\'on 2022)). Dr. Dionisio Y\'a\~nez has been supported through project CIAICO/2021/227 (Proyecto financiado por la Conselleria de Innovaci\'on, Universidades, Ciencia y Sociedad digital de la Generalitat Valenciana) and by grant PID2020-117211GB-I00 funded by MCIN/AEI/10.13039/501100011033.}

\date{}
\author{Sergio Amat,
\thanks{
 Departamento de Matem\'atica Aplicada y Estad\'{\i}stica.
  Universidad Polit\'ecnica de Cartagena (Spain).
e-mail:{\tt sergio.amat@upct.es}}\and
David Levin,
\thanks{
School of Mathematical Sciences. Tel-Aviv University, Tel-Aviv (Israel).
 e-mail:{\tt levindd@gmail.com}
}
\and Juan Ruiz-\'Alvarez, \thanks{ Corresponding Author. Departamento de Matem\'atica Aplicada y Estad\'{\i}stica.
  Universidad Polit\'ecnica de Cartagena (Spain).
e-mail:{\tt juan.ruiz@upct.es}}
\and Dionisio F. Y\'a\~nez \thanks{Departamento de Matem\'aticas. Universidad de Valencia, Valencia (Spain). e-mail:{\tt Dionisio.Yanez@uv.es}}
}
\maketitle

%\author[UPCT]{Sergio Amat}
%\ead{sergio.amat@upct.es}
%\author[TAU]{David Levin}
%\ead{levindd@gmail.com}
%\author[UPCT]{Juan Ruiz-\'Alvarez}
%\ead{juan.ruiz@upct.es}
%\author[UV]{Dionisio F. Y\'a\~nez}
%\ead{Dionisio.Yanez@uv.es}
%\address[UPCT]{Departamento de Matem\'atica Aplicada y Estad\'istica. Universidad  Polit\'ecnica de Cartagena. Cartagena (Spain).}
%\address[TAU]{School of Mathematical Sciences. Tel-Aviv University, Tel-Aviv (Israel).}
%\address[UV]{Departamento de Matem\'aticas. Universidad de Valencia, Valencia (Spain).}
%
%The correct dates will be entered by the editor.

\begin{abstract}
In this work, we present a new WENO B-spline-based quasi-interpolation algorithm. The novelty of this construction resides in the application of the WENO weights to the B-spline functions, that form a partition of unity, instead of the coefficients that multiply the B-spline functions of the spline. The result obtained conserves the smoothness of the original spline and presents adaption to discontinuities in the function. Another new idea that we introduce in this work is the use of different base weight functions from those proposed in classical WENO algorithms. Apart from introducing the construction of the new algorithms, we present theoretical results regarding the order of accuracy obtained at smooth zones and close to the discontinuity, as well as theoretical considerations about how to design the new weight functions. Through a tensor product strategy, we extend our results to several dimensions. In order to check the theoretical results obtained, we present an extensive battery of numerical experiments in one, two, and three dimensions that support our conclusions.
\end{abstract}
%\begin{keyword}
%WENO \sep high accuracy interpolation \sep improved adaption to discontinuities \sep generalization   \sep 41A05 \sep  41A10 \sep 65D05 \sep 65M06 \sep 65N06
%\end{keyword}
%\end{frontmatter}
%All acknowledgements should be placed in the back of the paper after Conclusions..
\vspace{-0.5cm}

\section{Introduction}
The WENO (weighted essentially non-oscillatory) method is an approximation strategy used for dealing with jump discontinuities that may occur in the solutions of hyperbolic PDEs.
The main idea of the WENO method is to build approximations using information from only one side of the discontinuity. This is achieved by utilizing certain smoothness indicators in order to construct a weighted approximation that uses data from smooth zones mainly. We note that in the development of numerical schemes for hyperbolic PDEs, the smoothness of the local approximation is less important and is less addressed.

A very brief and incomplete historical revision about the WENO algorithm can be found in what follows. Liu, Osher and Chan introduced for the first time the WENO algorithm in \cite{Liu}, and it was proposed as a local nonlinear convex combination of different interpolants, covering adjacent stencils. The values of the weights of this combination were functions of some smoothness indicators, calculated over the same stencil of the interpolants. These weights are designed so that the stencils affected by a discontinuity have a negligible contribution to the convex combination. Using a measure inspired by the total variation, Jiang and Shu introduced more suitable and efficient smoothness indicators in \cite{JiangShu}. The final objective of this technique was to provide high order of accuracy at smooth zones of the function while providing adaption to singularities using the same information as ENO (essentially non-oscillatory), \cite{MR881365, Harten1987231, SHU1988439, Shu1999, SHU198932, 
Arandiga2008225, TODOS, SM,CDM,Amat2007}. In the literature can be found many references about WENO algorithm. An incomplete list that the reader can visit to obtain more information about the state of the art of this algorithm can be: \cite{Shu1999,doi:10.1137/100791579, AMB, Henrick2005542, Castro20111766,RATHAN20181531,HUANG2018498,doi:10.1137/16M1087291,VANLITH2017529,Abgrall2016xxi,SUN201681,ZHU2016110,ZHANG20122245,Gerolymos:2009:VWS:1630172.1630406,doi:10.1002/fld.2168,doi:10.1137/100791002,doi:10.1002/num.21781,Arandiga:2014:WDM:2659923.2659932,li_liu_zhang_2015,Balsara:2016:ECW:3043104.3043133}, with a special emphasis in
\cite{Shu1998, doi:10.1137/070679065} and the references therein.

More recently, some investigation has been done about how to design WENO algorithms for splines \cite{AML,{ADLU}}. In this paper, we address the problem of approximating functions with jump discontinuities using their data on a uniform  grid.
To do so, we develop a new class of WENO-like methods which are based upon B-spline quasi-interpolation operators. There are several new ideas as well as new features presented within this development. The main idea is to apply the smoothness-based weights to the B-spline basis functions rather than to their coefficients. The second novelty is the introduction of a new weighting strategy challenging the classical one. The approximation method is defined with splines of any degree, and the resulting approximations are as smooth as the B-splines involved.

The method is first presented for the univariate case, whilst the extension to higher dimensions is straightforward using the tensor product strategy. An analysis  of the approximation error is presented, and extensive numerical tests and experiments are performed. Numerical results are also available for the approximation of functions with jump discontinuities in 2D and in 3D.

\section{Preliminaries}

In this section, we present the idea to construct a new non-linear operator which avoids the Gibbs-type phenomenon in the adjacent intervals close to a singularity and conserves the order of accuracy in the smooth zones. Using the same notation as \cite{explicitamatetal}, we consider a function $f$, a constant $h>0$, and the equidistant  data $\{f_n:=f(nh)\}$. For the approximation by a spline of degree $p$ we denote
$$f_{n,p}=(f_{n-\left\lfloor \frac{p}{2} \right\rfloor},\hdots,f_{n+\left\lfloor \frac{p}{2} \right\rfloor}),$$
and define the quasi-interpolation operator 
\begin{equation}\label{operatorQ}
Q_p(f)(x)=\sum_{n\in \mathbb{Z}}L_{p}(f_{n-\left\lfloor \frac{p}{2} \right\rfloor},\hdots,f_{n+\left\lfloor \frac{p}{2} \right\rfloor})B_p\left(\frac{x}{h}-n \right),
\end{equation}
where the different components are:
\begin{itemize} 
\item The $p$-degree B-spline, $B_p$, supported on
$\mathcal{I}_p =\left[-\frac{p+1}{2},\frac{p+1}{2}\right],$ with equidistant knots
$S_p = \left\{-\frac{p+1}{2},\hdots,\frac{p+1}{2}\right\}.$
\item The linear operator $L_p$, that is defined using the ideas presented in \cite{speleers} as:
\begin{equation}\label{operadorL}
L_p(f_{n,p})=\sum_{j=-\left\lfloor \frac{p}{2} \right\rfloor}^{\left\lfloor \frac{p}{2} \right\rfloor} c_{p,j} f_{n+j},
\end{equation}
with the coefficients $c_{p,j}$, $j=-\left\lfloor \frac{p}{2} \right\rfloor,\hdots,\left\lfloor \frac{p}{2} \right\rfloor$:
\begin{equation}\label{coeficientesc}
c_{p,j}=\sum_{l=0}^{\left\lceil \frac{p+1}{2} \right\rceil-1}\frac{t(2l+p+1,p+1)}{{{2l+p+1}\choose{p+1}}}\sum_{i=0}^{2l}\frac{(-1)^i}{i!(2l-i)!}\delta_{l-i+\left\lceil\frac{p+1}{2}\right\rceil,j+1+\left\lfloor\frac{p}{2}\right\rfloor},
\end{equation}
 where $\delta_{i,j}$ is the Kronecker delta function and $t(i,j)$ are the central factorial numbers of the first kind (see \cite{speleers} and \cite{butzeretal}), which can be computed recursively as:
\begin{equation*}
t(i,j)=\left\{
         \begin{array}{ll}
           0, & \hbox{if}\,\, j>i, \\
           1, & \hbox{if}\,\, j=i, \\
           t(i-2,j-2)-\left(\frac{i-2}{2}\right)^2t(i-2,j), & \hbox{if} \,\, 2\leq j <i,
         \end{array}
       \right.
\end{equation*}
with
$$t(i,0)=0,\quad t(i,1)=\prod_{l=1}^{i-1}\left(\frac{i}{2}-l\right),\,\, i\geq 2,$$
and $t(0,0)=t(1,1)=1, t(0,1)=t(1,0)=0$.
In Table \ref{tablaCs} we show some values of $c_{p,j}$.
\end{itemize}
%we call $B_p(x)$ the $p$-degree B-spline  Let $h>0$ be a constant and $\{f_n:=f(nh)\}$ the evaluation of $f$ in these nodes. We define the vector
%$$f_{n,p}=(f_{n-\left\lfloor \frac{p}{2} \right\rfloor},\hdots,f_{n+\left\lfloor \frac{p}{2} \right\rfloor}),$$
%where $\lfloor \cdot \rfloor$  is the floor function. The classical quasi-interpolation operator $Q_p$ is defined as
%\begin{equation}\label{operatorQ}
%Q_p(f)(x)=\sum_{n\in \mathbb{Z}}L_{p}(f_{n,p})B_p\left(\frac{x}{h}-n \right),
%\end{equation}
%with the linear operator $L_p$ defined using the ideas presented in \cite{speleers} as:

\begin{table}[htbp]
  \begin{center}
   \begin{tabular}{crrrrr}\hline
                             & $c_{1,j}$ & $c_{2,j}$ & $c_{3,j}$ & $c_{4,j}$& $c_{5,j}$    \\    \hline
                 $ j=0$      &  $1$  & $5/4 $    & $4/3$    & $319/192$  & $73/40$             \\
                 $ j=1$      &       & $-1/8 $   & $-1/6$   & $-107/288$ & $-7/15$               \\
                 $ j=2$      &       &           &          & $47/1152$  & $13/240$  \\\hline
   \end{tabular}
    \caption{$c_{p,j}$ values for $p=1,\hdots,5$. They are symmetric, i.e. $c_{p,j}=c_{p,-j}$, $j=0,\hdots,\left\lfloor\frac{p}{2}\right\rfloor$.}
    \label{tablaCs}
  \end{center}
\end{table}

The main power of the quasi-interpolator $Q_p$ lies in its capacity of obtaining explicit smooth, high-order, approximations to the function, as described in the following propositions (see for example \cite{speleers}):

\begin{proposition}\label{quasiintprop}
Consider a $p$ degree B-spline, $B_{p}$, with knots $S_p$, then the local operator $Q_{p}$ defined in Equation \eqref{operatorQ} is $p-1$ continuous.
\end{proposition}

\begin{proposition}\label{quasiintprop}
Consider a $p$ degree B-spline, $B_{p}$, with knots $S_p$, then the local operator $Q_{p}$ defined in Equation \eqref{operatorQ} reproduces $\Pi_{p}(\mathbb{R})$ and achieves $O(h^{p+1})$ approximation order to $C^{p+1}$ functions.
\end{proposition}

However, when the data present singularities or strong gradients, the reconstruction presents some artifacts close to the discontinuities (see e.g. \cite{AML,{ADLU}}). This Gibbs-type phenomenon can be avoided by constructing a new non-linear operator, as described in the next section.

\section{Non-linear WENO-based B-spline method}
This section aims to obtain a non-linear approximation operator by applying the arguments introduced for the WENO method to the operator $Q_p$. If we analyze Eq. \eqref{operatorQ} for a fixed point $x^*\in\mathbb{R}$, due to the compact support of the $B_p$ function, we can rewrite it as
\begin{equation}\label{operatorQ2}
Q_p(f)(x^*)=\sum_{n\in \mathbb{Z}}L_{p}(f_{n,p})B_p\left(\frac{x^*}{h}-n \right)=\sum_{k\in \mathcal{J}(x^*)}C^p_k(x^*)L_{p}(f_{k,p}),
\end{equation}
where
\begin{equation}\label{paraCs}
\mathcal{J}(x^*)=\left\{n\in\mathbb{Z}:B_p\left(\frac{x^*}{h}-n\right)> 0\right\}, \,\, \text{with} \,\, \left|\mathcal{J}(x^*)\right|\leq p+1, \quad C^p_k(x^*)=B_p\left(\frac{x^*}{h}-k\right)>0.
\end{equation}
We call the values $C^p_k(x^*)$ the B-optimal weights for $x^*$. Note that Eq. \eqref{operatorQ2} is a convex combination since
$\sum_{k\in \mathcal{J}(x^*)}C^p_k(x^*)=1.$
 We now substitute these linear weights with non-linear ones following the WENO strategy (see, e.g. \cite{JiangShu, Liu}). Thus,
\begin{equation}\label{operatorQ2omega}
Q^{\text{NL}}_p(f)(x^*)=\sum_{k\in \mathcal{J}(x^*)}\omega^p_k(x^*)L_{p}(f_{k,p}),
\end{equation}
being
\begin{equation}
\begin{aligned}
\omega^p_{k}(x^*)=\frac{\alpha_{k}^p(x^*)}{\sum_{j\in\mathcal{J}(x^*)} \alpha^p_j(x^*)},\quad \text{with} \quad \alpha_{k}^p={C}^p_{k}(x^*)\Psi({I}_{k,p}), \quad k\in \mathcal{J}(x^*),
\end{aligned}
\end{equation}
where $\Psi$ is a function that satisfies certain conditions to ensure maximum order at the smooth parts of the data. We will introduce these conditions within the next results.
In the WENO algorithm context, the function $\Psi$ has been typically defined as:
 $$ \Psi({I}_{k,p})=(\epsilon_h+{I}_{k,p})^{-t},$$
where $\epsilon_h$  is a parameter to avoid zero values in the denominator, for example, $\epsilon_h=h^2$, and the parameter $t$ is an integer that assures maximum order of accuracy close to the discontinuities.
Finally, the values ${I}_{k,p}$ are the smoothness indicators (see, e.g. \cite{JiangShu, Liu}) which are designed to detect which stencils are contaminated by a discontinuity. In our case, we choose the following:
\begin{equation}\label{sm:indicators}
{I}_{k,p}(f)=
\begin{cases}
(\Delta^p_k f)^2=\left(\sum_{j=-\frac{p}{2}}^{\frac{p}{2}}(-1)^{j+1} \binom{p}{j+\frac{p}{2}}f_{k+j}\right)^2, & \,\,\text{if $p$ is an even number},\\
(\Delta^{p-1}_k f)^2=\left(\sum_{j=-\left\lfloor \frac{p}{2}\right\rfloor}^{\left\lfloor \frac{p}{2}\right\rfloor}(-1)^{j+1} \binom{p-1}{j+\left\lfloor \frac{p}{2}\right\rfloor}f_{k+j}\right)^2, &  \,\,\text{if $p$ is an odd number}.\\
\end{cases}
\end{equation}
The properties of these operators are well-known, and they are essential to get optimal order of approximation at the smooth parts of the function, and to avoid the Gibbs phenomenon near singularities. We summarize the principal characteristics in the next proposition.
\begin{proposition}
Let $h>0$ and $I_{k,p}(f)$, $x^*\in\mathbb{R}$ the smoothness indicators defined in Eq. \eqref{sm:indicators} with $k\in\mathcal{J}(x^*)$ then:
\begin{itemize}
\item If $f\notin \mathcal{C}[h(k-\left\lfloor \frac{p}{2} \right\rfloor),h(k+\left\lfloor \frac{p}{2} \right\rfloor)]$ then
\begin{equation}\label{ordensuavidad}
I_{k,p}(f)=O(1),\quad h\to 0.
\end{equation}
\item If $f$ is a piecewise continuous function in $[h(k-\left\lfloor \frac{p}{2} \right\rfloor),h(k+\left\lfloor \frac{p}{2} \right\rfloor)]$ then
\begin{equation}\label{indicessuaves}
I_{k,p}(f)=O(h),\quad h\to 0.
%I_{k,p}(f)=\begin{cases}
%O(h^{2\min\{l,p\}}),& \text{if $p$ is an even number},\\
%O(h^{2\min\{l,p-1\}}),& \text{if $p$ is an odd number},\\
%\end{cases}
\end{equation}
%\item If $f\in \mathcal{C}^{l}[h(k-\left\lfloor \frac{p}{2} \right\rfloor),h(k+\left\lfloor \frac{p}{2} \right\rfloor)]$ with $l\geq p$ then
%$$I_{k,p}=\begin{cases}
%O(h^{2p}),& \text{if $p$ is an even number},\\
%O(h^{2p-2}),& \text{if $p$ is an odd number}.\\
%\end{cases}
%$$
%as $h\to 0$.
\item If $f\in \mathcal{C}^{p+1}[h(\min\mathcal{J}(x^*)-\left\lfloor \frac{p}{2} \right\rfloor),h(\max\mathcal{J}(x^*)+\left\lfloor \frac{p}{2} \right\rfloor)]$ then
\begin{equation}\label{distanciaentreindices}
I_{k,p}(f)-I_{l,p}(f)=\begin{cases}
O(h^{2p+1}),& \text{if $p$ is an even number},\\
O(h^{2p-1}),& \text{if $p$ is an odd number},\\
\end{cases}
\quad \forall\,\, k,l\in \mathcal{J}(x^*),
\end{equation}
as $h\to 0$.
\end{itemize}
\end{proposition}
\begin{demo}
The proof is straightforward using Taylor's expansion.
% If $0\leq l<p$ and $f\in\mathcal{C}^l[h(k-\left\lfloor \frac{p}{2} \right\rfloor),h(k+\left\lfloor \frac{p}{2} \right\rfloor)]$ then there exists $\xi_j \in [h(k-\left\lfloor \frac{p}{2} \right\rfloor),h(k+\left\lfloor \frac{p}{2} \right\rfloor)]$ such that:
%$$f_{k+j}=\sum_{i=0}^{l-1}\frac{f^{i)}(kh)}{i!}(jh)^i + \frac{f^{l)}(\xi_j)}{l!}(hj)^{l}$$
\end{demo}
%Note, that we only consider that the function is continuous $\mathcal{C}^{p+1}$ or it presents a jump discontinuity or a kink. 
For convenience, in the rest of the paper, we denote $I_{k,p}=I_{k,p}(f)$.

In order to guarantee that the behavior of the reconstruction obtained using the new operator is similar to the linear one at the smooth zones, we impose some conditions on the function $\Psi$. In the following proposition, we present the definition of the new non-linear approximation operator, together with the first requirement on the function $\Psi$.

\begin{proposition}\label{propo0}
Let $h>0$, $x^*\in\mathbb{R}$ and let $I_{k,p}$ be the smoothness indicators of $f$, with $k\in \mathcal{J}(x^*)$. Let $\Psi$ be a function such that \begin{equation}\label{condicionparapsi}
\Psi(I_{k,p})=\Psi(I_{l,p})(1+K_{k,l}h^{p+1}), \quad \forall k,l\in \mathcal{J}(x^*),
\end{equation}
with $K_{k,l}\in\mathbb{R}$ independent of $h$. Then,
\begin{equation}\label{pesos1}
\omega^p_{k}(x^*)=C^p_{k}(x^*)(1+O(h^{p+1})), \,\, \text{as $h\to 0$}.
\end{equation}
\end{proposition}
\begin{demo}
We fix a value $l\in \mathcal{J}(x^*)$, we have that for any $k\in \mathcal{J}(x^*)$ there exist $K_{k,l}\in\mathbb{R}$ such that
$$\alpha^p_k(x^*)={C}^p_{k}(x^*)\Psi(I_{k,p})={C}^p_{k}(x^*)\Psi(I_{l,p})(1+K_{k,l}h^{p+1})$$ %+{C}^p_{k}(x^*) K_{k,l} h^{p+1},$$
then from $\sum_{j\in\mathcal{J}(x^*)} {C}^p_{j}(x^*)=1$ we get
$$\sum_{j\in\mathcal{J}(x^*)} \alpha^p_j(x^*)=\sum_{j\in\mathcal{J}(x^*)} {C}^j_{p}(x^*)\Psi({I}_{j,p})=\sum_{j\in\mathcal{J}(x^*)} {C}^p_{j}(x^*)\Psi(I_{l,p})(1+K_{j,l}h^{p+1})=\Psi(I_{l,p})(1+\tilde{K}h^{p+1}).$$
Finally,
$$\omega^p_k(x^*)=\frac{\alpha_{k}^p(x^*)}{\sum_{k\in\mathcal{J}(x^*)} \alpha^p_k(x^*)}=\frac{{C}^p_{k}(x^*)\Psi(I_{l,p})(1+K_{k,l}h^{p+1})}{\Psi(I_{l,p})(1+\tilde{K}h^{p+1})} ={C}^p_{k}(x^*)(1+O(h^{p+1})). $$
\end{demo}

With these ingredients, we have the following main result concerning the approximation order of the non-linear operator at the smooth regions of $f$.
\begin{proposition}\label{propoQ}
Let be $h>0$, $x^*\in \mathbb{R}$, $Q^{\text{NL}}_p$ the operator defined in Eq. \eqref{operatorQ2omega} and $\Psi$ a function which satisfies Eq. \eqref{condicionparapsi}. Then, if $f\in\mathcal{C}^{p+1}[h(k-\left\lfloor \frac{p}{2} \right\rfloor),h(k+\left\lfloor \frac{p}{2} \right\rfloor)]$, for all $k\in\mathcal{J}(x^*)$ then:
$$|Q^{\text{NL}}_p(f)(x^*)-f(x^*)|=O(h^{p+1}),\,\,\text{as $h\to 0$}.$$

\end{proposition}
\begin{demo}
The proof is straightforward, by Prop. \ref{propo0}:
\begin{equation}
\begin{split}
|Q^{\text{NL}}_p(f)(x^*)-Q_p(f)(x^*)|&=\left|\sum_{k\in \mathcal{J}(x^*)}C^p_k(x^*)L_{p}(f_{k,p})-\sum_{k\in \mathcal{J}(x^*)}\omega^p_k(x^*)L_{p}(f_{k,p})\right|\\
&=\left|\sum_{k\in \mathcal{J}(x^*)}(C^p_k(x^*)-\omega^p_k(x^*))L_{p}(f_{k,p})\right|\\
&=\left|\sum_{k\in \mathcal{J}(x^*)}O(h^{p+1})L_{p}(f_{k,p})\right|=O(h^{p+1}),\\
\end{split}
\end{equation}
and from $f\in\mathcal{C}^{p+1}[h(k-\left\lfloor \frac{p}{2} \right\rfloor),h(k+\left\lfloor \frac{p}{2} \right\rfloor)]$, for all $k\in\mathcal{J}(x^*)$, using Proposition \ref{quasiintprop}, we get:
$$|Q^{\text{NL}}_p(f)(x^*)-f(x^*)|\leq |Q^{\text{NL}}_p(f)(x^*)-Q_p(f)(x^*)|+|Q_p(f)(x^*)-f(x^*)|=O(h^{p+1}).$$
\end{demo}

With the condition given in Proposition \ref{propoQ}  the accuracy should be optimal, at least, when the data do not present singularities. In the following lemma, we present
sufficient conditions for constructing a function $\Psi$ with the required properties.
\begin{lemma}\label{lemma1}
Let $h>0$, $x^*\in\mathbb{R}$, $\psi \in \mathcal{C}^2$ with $\psi(x)\neq 0$ for all $x\in \mathbb{R}_+$, and let $I_{k,p}$ be smoothness indicators of $f$, with $f\in \mathcal{C}^{p+1}[h(\min\mathcal{J}(x^*)-\left\lfloor \frac{p}{2} \right\rfloor),h(\max\mathcal{J}(x^*)+\left\lfloor \frac{p}{2} \right\rfloor)]$, such that
\begin{align*}
&\psi(I_{k,p})=O(h^s),\,\,  \forall \,\, k\in \mathcal{J}(x^*), \\
&\psi'(x)=O(h^{q_1}), \,\, \forall \,\, x\in [\min(I_{l,p},I_{k,p}),\max(I_{l,p},I_{k,p})],\\
&\psi''(x)=O(h^{q_2}), \,\, \forall \,\, x\in [\min(I_{l,p},I_{k,p}),\max(I_{l,p},I_{k,p})],\\
&I_{l,p}-I_{k,p}=O(h^{s_1}), \quad \forall k,l\in \mathcal{J}(x^*),
\end{align*}
as $h\to 0$, with $s_1,s\in\mathbb{N}$, $s_1\geq s$ and $q_1,q_2\in\mathbb{Z}$, $s_1+q_2\geq q_1$. Then, the function $\Psi(x)=\frac{1}{\psi(x)}$ satisfies:
$$\Psi(I_{k,p})=\Psi(I_{l,p})(1+O(h^{q_1+s_1-s})),\,\, \text{as $h\to 0$}.$$

\end{lemma}

\begin{demo}
As $\psi \in \mathcal{C}^2$ then $\Psi \in \mathcal{C}^2$ and, by Taylor's expansions, there exists $\xi\in (\min(I_{l,p},I_{k,p}),\max(I_{l,p},I_{k,p}))$ such that
\begin{equation*}
\begin{split}
\Psi(I_{k,p})&=\Psi(I_{l,p})+\Psi'(I_{l,p})(I_{k,p}-I_{l,p})+\frac{\Psi''(\xi)}{2}(I_{k,p}-I_{l,p})^2\\
&=\Psi(I_{l,p})+\frac{\psi'(I_{l,p})}{\psi(I_{l,p})^2}(I_{l,p}-I_{k,p})+\frac{2
   \psi'(\xi)^2-\psi(\xi) \psi''(\xi)}{2 \psi(\xi)^3}(I_{l,p}-I_{k,p})^2\\
&=\Psi(I_{l,p})\left(1+\frac{\psi'(I_{l,p})}{\psi(I_{l,p})}(I_{l,p}-I_{k,p})+\frac{2
   \psi'(\xi)^2-\psi(\xi) \psi''(\xi)}{2 \psi(\xi)^3}\psi(I_{l,p})(I_{l,p}-I_{k,p})^2\right)\\
&=\Psi(I_{l,p})\left(1+\frac{O(h^{q_1})}{O(h^{s})}O(h^{s_1})+\frac{2
   O(h^{2q_1})-O(h^{s}) O(h^{q_2})}{2 O(h^{3s})}O(h^s)O(h^{2s_1})\right)\\
&=\Psi(I_{l,p})\left(1+O(h^{\min(q_1+s_1-s, 2(q_1+s_1-s),2s_1+q_2-s)})\right)\\
&=\Psi(I_{l,p})\left(1+O(h^{q_1+s_1-s})\right).\\
\end{split}
\end{equation*}%
%$$ \Psi(I_{l,p})\left(\frac{(I_{l,p}-I_{k,p}) \psi'(I_{l,p})}{\psi(I_{l,p})}+\frac{(I_{l,p}-I_{k,p})^2 \left(2
%   \psi'(I_{l,p})^2-\psi(I_{l,p}) \psi''(I_{l,p})\right)}{2 \psi(I_{l,p})^2}\right)+O\left((I_{k,p}-I_{l,p})^3\right)$$
\end{demo}

Note that, if we determine $\Psi(x)=C\neq 0$, for all $x\in\mathbb{R}$ we get good behaviour when the data do not present any singularity. This is clear because, in this case, $C_k^p=\omega_k^p$, for all $k$ and, hence, $Q_p=Q_p^{\text{NL}}$. Therefore, we impose different conditions to adapt our new operator to the discontinuities in the following result.

\begin{proposition}\label{propo2}
Let us consider $h>0$, $\psi \in \mathcal{C}^2$ with $\psi(x)\neq 0$ for all $x\in \mathbb{R}$, $\Psi(x)=\frac{1}{\psi(x)}$, and let $I_{n,p}$ be smoothness indicators of $f$. If the function $f$ presents a jump or a kink (a jump in the first order derivative) in at least one interval, $n\in\mathcal{J}(x^*)$, and, as $h\to 0$,
\begin{equation}\label{equationpsi1}
\begin{split}
&\psi(I_{k,p})=O(h^s),\quad  \text{if} \,\, f\in\mathcal{C}^{p+1}[h(k-\left\lfloor \frac{p}{2} \right\rfloor),h(k+\left\lfloor \frac{p}{2} \right\rfloor)], \\
&\psi(I_{k,p})=O(h^q),\quad  \text{otherwise},\\
\end{split}
\end{equation}
with $s\in\mathbb{N}$, $q\in\mathbb{Z}$, $s> q$,  then:
$$\omega_k^p(x^*)=
\begin{cases}
O(1),&  \text{if} \,\, f\in\mathcal{C}^{p+1}[h(k-\left\lfloor \frac{p}{2} \right\rfloor),h(k+\left\lfloor \frac{p}{2} \right\rfloor)], \\
O(h^{s-q}),&  \text{otherwise}. \\
\end{cases}$$
\end{proposition}
\begin{demo}
By \eqref{equationpsi1}, we have that
$$\alpha_k^p(x^*)=
\begin{cases}
\frac{C^p_k(x^*)}{\psi(I_{k,p})}=O(h^{-s}),&  f\in\mathcal{C}^{p+1}[h(k-\left\lfloor \frac{p}{2} \right\rfloor),h(k+\left\lfloor \frac{p}{2} \right\rfloor)], \\
\frac{C^p_k(x^*)}{\psi(I_{k,p})}=O(h^{-q}),&  \text{otherwise.}
\end{cases}$$
Thus, as $s> q$ we calculate $\sum_{k\in \mathcal{J}(x^*)} \alpha^p_k(x^*)=O(h^{\min\{-s,-q\}})=O(h^{-s})$. If $f\in\mathcal{C}^{p+1}[h(k-\left\lfloor \frac{p}{2} \right\rfloor),h(k+\left\lfloor \frac{p}{2} \right\rfloor)]$ we have that
$$\omega^p_k(x^*)=\frac{\alpha_k^p(x^*)}{\sum_{k\in \mathcal{J}(x^*)} \alpha^p_k(x^*)}=\frac{O(h^{-s})}{O(h^{-s})}=O(1),$$
and if $f$ presents a jump or a kink in $[h(k-\left\lfloor \frac{p}{2} \right\rfloor),h(k+\left\lfloor \frac{p}{2} \right\rfloor)]$ we get
$$\omega^p_k(x^*)=\frac{\alpha_k^p(x^*)}{\sum_{k\in \mathcal{J}(x^*)} \alpha^p_k(x^*)}=\frac{O(h^{-q})}{O(h^{-s})}=O(h^{s-q}).$$
\end{demo}

Finally, we prove that under these conditions we have approximation order $O(h)$ close to the singularities.

\begin{proposition}\label{ordenh}
Let us consider $h>0$, $n_0\in\mathbb{N}$, $p\geq 2$, $\psi \in \mathcal{C}^2$ satisfying the hypothesis of Proposition \ref{propo2}, and $x^*\in [n_0 h,(n_0+1)h]$, $\xi \in ((n_0+1)h,(n_0+p)h)$ if $p$ is an odd number and $\xi \in ((n_0+1)h,(n_0+p+1)h)$ otherwise. We consider a discontinuous function, $f$ given by
$$f(x)=\begin{cases}
f^+(x), & x< \xi,\\
f^-(x), & x\geq \xi,\\
\end{cases}
$$
with $f^+\in \mathcal{C}^{p+1}]-\infty,\xi]$, $f^-\in \mathcal{C}^{p+1}[\xi,+\infty)$, with $f^+(\xi)> f^-(\xi)$ then
$$|Q^{\text{NL}}_p(f)(x^*)-f(x^*)|=O(h).$$
\end{proposition}
\begin{demo}
We suppose that $p$ is an odd number, the proof is similar for an even $p$. Let $l\in\mathbb{N}$ with $1\leq l\leq p$ such that $\xi \in ((n_0+l)h,(n_0+l+1)h)$ and
$k\in \left\{n_0-\left\lfloor \frac{p}{2} \right\rfloor,\hdots,n_0-\left\lfloor \frac{p}{2} \right\rfloor+l\right\}$. If $f$ is smooth in $[h(k-\left\lfloor \frac{p}{2} \right\rfloor),h(k+\left\lfloor \frac{p}{2} \right\rfloor)]$, then
$$f_{k+j}=f((k+j)h)=f(kh)+O(h)=f_k+O(h)=f(x^*)+O(h), \quad -\left\lfloor \frac{p}{2} \right\rfloor\leq j \leq \left\lfloor \frac{p}{2} \right\rfloor.$$
From the definition of the operator $L_p$, Eq. \eqref{operadorL}, we have that
\begin{equation}\label{eqlemapauso}
L_p(f_{k,p})=\sum_{j=-\left\lfloor \frac{p}{2} \right\rfloor}^{\left\lfloor \frac{p}{2} \right\rfloor} c_{p,j} f_{k+j}=\sum_{j=-\left\lfloor \frac{p}{2} \right\rfloor}^{\left\lfloor \frac{p}{2} \right\rfloor} c_{p,j}(f(x^*)+O(h))=f(x^*)+O(h),
\end{equation}
since $\sum_{j=-\left\lfloor \frac{p}{2} \right\rfloor}^{\left\lfloor \frac{p}{2} \right\rfloor} c_{p,j}=1$. Now, if $k'\in \left\{n_0-\left\lfloor \frac{p}{2} \right\rfloor+l+1,\hdots,n_0+\frac{p+1}{2}\right\}$ we have that
$$((n_0+l)h,(n_0+l+1)h)\subseteq \left[\left(k'-\left\lfloor \frac{p}{2}\right\rfloor\right)h, \left(k'+\left\lfloor \frac{p}{2}\right\rfloor\right)h\right].$$
This expression holds since, if $k'=n_0-\left\lfloor \frac{p}{2} \right\rfloor+l+1$, we obtain
\begin{equation*}
k'-\left\lfloor \frac{p}{2}\right\rfloor=n_0+l-p, \quad k'+\left\lfloor \frac{p}{2}\right\rfloor=n_0+l+1,\\
\end{equation*}
and if $k'=n_0+\frac{p+1}{2}$ we get
\begin{equation*}
k'-\left\lfloor \frac{p}{2}\right\rfloor=n_0+1, \quad k'+\left\lfloor \frac{p}{2}\right\rfloor=n_0+p.\\
\end{equation*}
Therefore, by Proposition \ref{propo2}, if  $k'\in \left\{n_0-\left\lfloor \frac{p}{2} \right\rfloor+l+1,\hdots,n_0+\frac{p+1}{2}\right\}$ then $\omega_{k'}^p(x^*)=O(h^{s-q})$, $s-q\geq 1$. Also, as  $\mathcal{J}(x^*)=\{n_0-\left\lfloor \frac{p}{2} \right\rfloor,\hdots,n_0+\frac{p+1}{2}\}$, notice that
\begin{equation}\label{sumapesosparalema}
1=\sum_{k\in \mathcal{J}(x^*)}\omega^p_k(x^*)=\sum_{k=n_0-\left\lfloor \frac{p}{2} \right\rfloor}^{n_0+\frac{p+1}{2}}\omega^p_k(x^*)=\sum_{k=n_0-\left\lfloor \frac{p}{2} \right\rfloor}^{n_0-\left\lfloor \frac{p}{2} \right\rfloor+l}\omega^p_k(x^*) + \sum_{k=n_0-\left\lfloor \frac{p}{2} \right\rfloor+l+1}^{n_0+\frac{p+1}{2}}\omega^p_k(x^*)=\sum_{k=n_0-\left\lfloor \frac{p}{2} \right\rfloor}^{n_0-\left\lfloor \frac{p}{2} \right\rfloor+l}\omega^p_k(x^*) + O(h^{s-q}).\\
\end{equation}

Finally, by \eqref{eqlemapauso} and \eqref{sumapesosparalema}, we have that
\begin{equation}
\begin{split}
Q^{\text{NL}}_p(f)(x^*)&=\sum_{k\in \mathcal{J}(x^*)}\omega^p_k(x^*)L_{p}(f_{k,p})=\sum_{k=n_0-\left\lfloor \frac{p}{2} \right\rfloor}^{n_0+\frac{p+1}{2}}\omega^p_k(x^*)L_{p}(f_{k,p})\\
&=\sum_{k=n_0-\left\lfloor \frac{p}{2} \right\rfloor}^{n_0-\left\lfloor \frac{p}{2} \right\rfloor+l}\omega^p_k(x^*)L_{p}(f_{k,p}) + \sum_{k=n_0-\left\lfloor \frac{p}{2} \right\rfloor+l+1}^{n_0+\frac{p+1}{2}}\omega^p_k(x^*)L_{p}(f_{k,p})\\
&=\sum_{k=n_0-\left\lfloor \frac{p}{2} \right\rfloor}^{n_0-\left\lfloor \frac{p}{2} \right\rfloor+l}\omega^p_k(x^*)(f(x^*)+O(h)) + \sum_{k=n_0-\left\lfloor \frac{p}{2} \right\rfloor+l+1}^{n_0+\frac{p+1}{2}}O(h^{s-q})L_{p}(f_{k,p})\\
&=f(x^*)+O(h)+O(h^{s-q})=f(x^*)+O(h).
\end{split}
\end{equation}
\end{demo}

\subsection{The design of the function $\Psi$}\label{psi}

In this subsection, we present some possible options to obtain a function $\Psi$ which satisfies the above-mentioned conditions.
We know, by Eq. \eqref{distanciaentreindices}, that $s_1=2p+1$ if $p$ is even, $s_1=2p-1$ if it is odd, and that $\Psi=1/\psi$ with $\psi\in\mathcal{C}^2(\mathbb{R})$ and $\psi(x)\neq 0$, for all $x\in\mathbb{R}_+$. Therefore, we only summarize the conditions of the function $\psi$:
\begin{enumerate}[label=(C\arabic*)]
\item $\psi(I_{k,p})=O(h^s),\,\,  \forall \,\, k\in \mathcal{J}(x^*)$, if $f\in \mathcal{C}^{p+1}[h(k-\left\lfloor \frac{p}{2} \right\rfloor),h(k+\left\lfloor \frac{p}{2} \right\rfloor)]$.
\item $\psi'(x)=O(h^{q_1}), \,\, \forall \,\, x\in [\min(I_{l,p},I_{k,p}),\max(I_{l,p},I_{k,p})]$ if $f\in \mathcal{C}^{p+1}[h(\min\mathcal{J}(x^*)-\left\lfloor \frac{p}{2} \right\rfloor),h(\max\mathcal{J}(x^*)+\left\lfloor \frac{p}{2} \right\rfloor)]$.
\item $\psi''(x)=O(h^{q_2}), \,\, \forall \,\, x\in [\min(I_{l,p},I_{k,p}),\max(I_{l,p},I_{k,p})]$ if $f\in \mathcal{C}^{p+1}[h(\min\mathcal{J}(x^*)-\left\lfloor \frac{p}{2} \right\rfloor),h(\max\mathcal{J}(x^*)+\left\lfloor \frac{p}{2} \right\rfloor)]$.
\item $\psi(I_{k,p})=O(h^q)$ if $f\notin \mathcal{C}[h(k-\left\lfloor \frac{p}{2} \right\rfloor),h(k+\left\lfloor \frac{p}{2} \right\rfloor)]$.
\end{enumerate}
Note that the first three conditions are devoted to obtaining optimal order, so, it is necessary that:
$$s_1+q_1-s\geq p+1, \quad s_1\geq s, \quad s_1+q_2\geq q_1.$$
The fourth condition is to adapt the operator to deal with discontinuities. In this case, the constraint is $q<s$. We suppose that if the function is not continuous then the function presents a jump discontinuity in at least one stencil, i.e., $I_{k,p}=O(1)$ as $h\to 0$.

We start with the classical smoothness indicators designed by Jiang and Shu in \cite{JiangShu}, that  where analyzed and modified in \cite{ABM}. Thereby, we consider the function $\psi_s(x)=h^2+x$. In this case, we have $s=2$, $q_1=q_2=0$ and $q=0$. Therefore, all the conditions are satisfied. The behavior of the resulting approximation using the new operator is adequate. 

Our second example is the function $\psi_c(x)=C+x/h$. It is clear that
$s=1$, its derivative is $\psi'(x)=1/h$ then $q_1=-1$ and $q_2=0$. Finally, in the presence of a jump discontinuity, we have that $q=-1$.
Again, all the conditions are fulfilled. For this function, we should pre-determine a constant $C$. If it is very big, the efficiency of the algorithm in the smooth parts increases, however, the capacity of detection of discontinuities is reduced. If $C$ is chosen to be small, the absolute error in the smooth part is higher than the one obtained when $C$ is bigger. We choose $C=1$ and we perform some numerical experiments in the next section. 

We introduce the last example to avoid determining constants. We define the function $\psi_d(x)=e^{x/h}$ and calculate the constants:
$s=0$, $q_1=-1$, $q_2=-2$. If $f\notin \mathcal{C}[h(k-\left\lfloor \frac{p}{2} \right\rfloor),h(k+\left\lfloor \frac{p}{2} \right\rfloor)]$ then 
$\psi_d(I_{k,p})=e^{O(1/h)}\geq O(h^q),\,\,q\in\mathbb{Z}_{-}$. A numerical accuracy comparison with these different options of $\psi$ is presented in section \ref{expnum}, where some advantage is observed when using the function $\psi_d$. Also, $\psi_d$ seems to provide more sharpness in the approximations close to discontinuities. 

\section{Generalization to several dimensions}\label{dimensions}

One of the main advantages of the new method is its ease of generalization to higher dimensions. In this section, we present the general method in $k$-dimensions. For this, we use the notation presented in \cite{explicitamatetal} (also see \cite{speleers}).
Let $h>0$, we consider a real function $f:\mathbb{R}^k\to \mathbb{R}$, the vectors $\mathbf{p}=(p_1,\hdots,p_k)$, $\mathbf{n}=(n_1,\hdots,n_k)$ and  the values
$$\{f_{\mathbf{n}}=f({n_{1} h,\hdots,n_{k}h})\}.$$
We denote as:
$$f_{\mathbf{n},\mathbf{p}}=\left\{f_{(n_{1}+j_1,\hdots,n_{k}+j_k)}:0\leq |j_l|\leq \left\lfloor \frac{p_l}{2}\right\rfloor, l=1,\hdots,k\right\}.$$
It is clear that if $k=1$, $f_{n,p}=(f_{n-\left\lfloor \frac{p}{2} \right\rfloor},\hdots,f_{n+\left\lfloor \frac{p}{2} \right\rfloor})$. Thus, the quasi-linear
operator in $n$ dimensions for a point $\mathbf{x}^*=(x^*_1,\hdots,x^*_k)\in\mathbb{R}^k$ is:
\begin{equation}\label{operatorQmulti}
Q_{\mathbf{p}}(f)(\mathbf{x}^*)=\sum_{\mathbf{n}\in\mathbb{Z}^k} L_{\mathbf{p}}(f_{\mathbf{n},\mathbf{p}})B_{\mathbf{p}}\left(\frac{\mathbf{x}^*}{h}-\mathbf{n}\right),
\end{equation}
where $L_{\mathbf{p}}$ is the tensor product of $L_p$, Equation \eqref{operadorL}, i.e.
\begin{equation}\label{operadorLmulti}
L_{\mathbf{p}}(f_{\mathbf{n},\mathbf{p}})=\sum_{j_1=-\left\lfloor \frac{p_1}{2} \right\rfloor}^{\left\lfloor \frac{p_1}{2} \right\rfloor}\hdots\sum_{j_k=-\left\lfloor \frac{p_k}{2} \right\rfloor}^{\left\lfloor \frac{p_k}{2} \right\rfloor} c_{p_1,j_1}\hdots c_{p_k,j_k} f_{n_1+j_1,\hdots,n_k+j_k},
\end{equation}
with $c_{p,j}$ defined in Equation \eqref{coeficientesc}, and
$$B_{\mathbf{p}}\left(\frac{\mathbf{x}^*}{h}-\mathbf{n}\right)=\prod_{l=1}^kB_{p_{l}}\left(\frac{x_l^*}{h}-n_l\right).$$
Using the same notation as in Eq. \eqref{paraCs}, we have that:
$$B_{\mathbf{p}}\left(\frac{\mathbf{x}^*}{h}-\mathbf{n}\right)=\prod_{l=1}^kC^{p_{l}}_{n_l}\left({x_l^*}\right).$$
To generalize the operator $Q^{\text{NL}}_p$, Eq. \eqref{operatorQ2omega}, we replace each value $C^{p_{l}}_{n_l}\left({x_l^*}\right)$ by its corresponding $\omega^{p_{l}}_{n_l}\left(x^*_l\right)$ in the same way as we have developed above. If we denote:
$$\mathcal{W}_{\mathbf{p}}\left(\mathbf{x}^*\right)=\prod_{l=1}^k\omega^{p_{l}}_{n_l}\left({x_l^*}\right),$$
the new non-linear operator is:
\begin{equation}\label{operatorQmultinolineal}
Q^{\text{NL}}_{\mathbf{p}}(f)(\mathbf{x}^*)=\sum_{\mathbf{n}\in\mathbb{Z}^k} L_{\mathbf{p}}(f_{\mathbf{n},\mathbf{p}})\mathcal{W}_{\mathbf{p}}\left(\mathbf{x}^*\right).
\end{equation}
As we can see, a tensor product strategy has been used. Therefore, if the function is sufficiently smooth the order of accuracy is $1+\min_{1\leq l\leq k} \{p_l\}$, and $O(h)$ when the function presents some singularity.

\section{Numerical experiments}\label{expnum}
%Todos los experimentos en: /Users/j/Documents/articulos_a_medias/Levin/bsplines_WENO_pph/graficas_disc_articulo.m

In this section, we analyze numerically some of the properties of the algorithms proposed in previous sections. More specifically, we check the behavior of the proposed WENO B-spline-based algorithms at smooth zones, and close to the discontinuities. We are interested in studying the numerical order and the absence of the Gibbs-type phenomenon. We also present results for bivariate and trivariate data, obtained by applying a tensor product strategy, as explained in Section \ref{dimensions}.

\subsection{Accuracy at smooth zones}
In this section, we analyze the behavior of the new algorithm and the classical one at regions of smoothness in order to check the numerical accuracy reached. To do so, we consider the smooth function
\begin{equation}\label{pol}
f(x)=x^6+x^3-3x^2,\ x \in [0,1],
\end{equation}
and perform a grid refinement analysis. The error in the approximated solution is denoted by $E^l$  and it is computed using a grid size $h_l$ using the infinity norm inside the considered interval,
\begin{equation}\label{linf}
E^l=||f^l-\tilde{f^l}||_{L^\infty}.
\end{equation}
We can see that $l$ represents the step of the grid refinement analysis, which consists in reducing the grid size $h$ to $h/2$ when going from $l$ to $l+1$. 
The order of accuracy is obtained in general as,
\begin{equation}\label{orden}
O^l=\log_2\left(\frac{E^l}{ E^{l+1}}\right).
\end{equation}

The results for the quadratic spline are presented in table \ref{tabla_quadratic}, for the cubic in table \ref{tabla_cubic}, for the quartic in table \ref{tabla_quartic} and for the quintic in table \ref{tabla_quintic}. All the tables show the results for the functions $\psi_s, \psi_c, \psi_d$, presented in Subsection \ref{psi}, and the classical spline. We can see how all the versions of the algorithm attain the correct numerical accuracy when applied to a smooth function: order $O(h^3)$ for the quadratic spline, order $O(h^4)$ for the cubic, and so on.

\begin{table}[!ht]
\begin{center}
%Programas en: /Users/j/Documents/articulos_a_medias/Levin/bsplines_WENO_pph/graficas_disc_articulo.m
\resizebox{15cm}{!} {
\begin{tabular}{|c|c|c|c|c|c|c|c|c|c|c|c|c|}
\hline $m=2^l$ & $2^4$ &  $2^5$           & $2^6$          &  $2^7$       & $2^8$ &  $2^9$   & $2^{10}$ %& $2^{11}$ & $2^{12}$
            \\
\hline
Error WENO-$\psi_s$ quadratic spline ($E^l$)&   2.1439e-02 & 3.1585e-03 & 2.9446e-04 & 2.0690e-05 & 1.3538e-06&8.9398e-08 & 6.0681e-09
            \\
\hline
$O^l$ & - & 2.763 & 3.423 & 3.831 & 3.934 & 3.921& 3.881
        \\
\hline        
Error WENO-$\psi_c$ quadratic spline ($E^l$)&  6.2029e-03 & 2.0395e-04 & 9.5224e-06 & 6.6798e-07 & 6.5957e-08&7.6283e-09 & 9.3011e-10
            \\
\hline
$O^l$ &- & 4.927 & 4.421 & 3.833 & 3.34 & 3.112 & 3.036
        \\        
\hline
Error WENO-$\psi_d$ quadratic spline ($E^l$)&   8.0140e-03 & 2.0957e-04 & 9.5404e-06 & 6.6805e-07 & 6.5957e-08&7.6283e-09 & 9.3011e-10
            \\
\hline
$O^l$ & - & 5.257 & 4.457 & 3.836 & 3.34 & 3.112& 3.036
        \\
\hline
\hline
Error class. quadratic spline ($E^l$)&    4.1225e-04 & 3.9638e-05 & 4.3232e-06 & 5.0409e-07 & 6.0818e-08& 7.4674e-09& 9.2508e-10
            \\
\hline
$O^l$ & - & 3.379 & 3.197 &  3.1 & 3.051 & 3.026 &3.013
        \\
\hline
\end{tabular}
}
\caption{Grid refinement analysis for the accuracy of the quadratic WENO spline and the classical quadratic spline using the infinity norm. The original data has been sampled from the function in (\ref{pol}). We show the results for the functions $\psi_s, \psi_c, \psi_d$, presented in Subsection \ref{psi}, and the classical spline. The low-resolution nodes have been sampled with $m=2^l$ nodes and the high-resolution data with $11(m-1)+m$ points.}\label{tabla_quadratic}
\end{center}
\end{table}

\begin{table}[!ht]
\begin{center}
%Programas en: /Users/j/Documents/articulos_a_medias/Levin/bsplines_WENO_pph/graficas_disc_articulo.m
\resizebox{15cm}{!} {
\begin{tabular}{|c|c|c|c|c|c|c|c|c|c|c|c|c|}
\hline $m=2^l$ & $2^4$ &  $2^5$           & $2^6$          &  $2^7$       & $2^8$ &  $2^9$   & $2^{10}$ %& $2^{11}$ & $2^{12}$
            \\
\hline
Error WENO-$\psi_s$ cubic spline ($E^l$)&       5.7275e-03  & 4.4580e-04  & 1.1572e-04 &  1.5515e-05&   1.3342e-06  & 9.6355e-08 &  6.4480e-09
            \\
\hline
$O^l$ & -&      3.6834  & 1.9457 &  2.8990  & 3.5396  & 3.7914  & 3.9014
        \\
\hline
Error WENO-$\psi_c$ cubic spline ($E^l$)&    3.8744e-04 & 9.4390e-06 & 2.6245e-06 & 1.6338e-07 & 7.6510e-09&3.4392e-10 & 1.6146e-11
            \\
\hline
$O^l$ & - & 5.359 & 1.847 & 4.006 & 4.416 & 4.476 & 4.413
        \\                
\hline
Error WENO-$\psi_d$ cubic spline ($E^l$)&    3.9653e-04 & 9.4891e-06 & 2.6288e-06 & 1.6342e-07 & 7.6513e-09&3.4392e-10 & 1.6147e-11
            \\
\hline
$O^l$ & - & 5.385 & 1.852 & 4.008 & 4.417 & 4.476 & 4.413
        \\
\hline
\hline
Error class. cubic spline ($E^l$)&    1.0716e-04  & 8.6705e-06  & 5.9865e-07  & 3.9087e-08  & 2.4935e-09 &  1.5740e-10  & 9.8859e-12
            \\
\hline
$O^l$ & -&  3.6275 &  3.8563 &  3.9370  & 3.9704 &  3.9857   &3.9929
        \\
\hline
\end{tabular}
}
\caption{Grid refinement analysis for the accuracy of the cubic WENO spline and the classical cubic spline using the infinity norm. The original data has been sampled from the function in (\ref{pol}). We show the results for the functions $\psi_s, \psi_c, \psi_d$, presented in Subsection \ref{psi}, and the classical spline. The low-resolution nodes have been sampled with $m=2^l$ nodes and the high-resolution data with $10(m-1)+m$ points.}\label{tabla_cubic}
\end{center}
\end{table}

\begin{table}[!ht]
\begin{center}
%Programas en: /Users/j/Documents/articulos_a_medias/Levin/bsplines_WENO_pph/graficas_disc_articulo.m
\resizebox{15cm}{!} {
\begin{tabular}{|c|c|c|c|c|c|c|c|c|c|c|c|c|}
\hline $m=2^l$ & $2^4$ &  $2^5$           & $2^6$          &  $2^7$       & $2^8$ &  $2^9$   & $2^{10}$ %& $2^{11}$ & $2^{12}$
            \\
\hline
Error WENO-$\psi_s$ quartic spline ($E^l$)&   2.5613e-04 & 7.5817e-07 & 2.5238e-09 & 1.1877e-11 & 1.4899e-13 & 4.8850e-15 & 1.1102e-15
            \\
\hline
$O^l$ & -&   8.4 & 8.231 & 7.731 & 6.317 & 4.931 & 2.138 
        \\
\hline
Error WENO-$\psi_c$ quartic spline ($E^l$)&  1.6809e-05 & 1.8213e-08 & 2.2147e-10 & 5.7438e-12 & 1.5721e-13 & 4.8850e-15 & 1.1102e-15
            \\
\hline
$O^l$ & -&  9.85 & 6.362 & 5.269 & 5.191 & 5.008 & 2.138
        \\        
\hline
Error WENO-$\psi_d$ quartic spline ($E^l$)&   1.6823e-05 & 1.8214e-08 & 2.2147e-10 & 5.7438e-12 & 1.5721e-13 & 4.8850e-15 & 1.1102e-15
            \\
\hline
$O^l$ & -&   9.851 & 6.362 & 5.269 & 5.191 & 5.008 & 2.138
        \\
\hline
\hline
Error class. quartic spline ($E^l$)&   7.0262e-07 & 1.1659e-08 & 2.4046e-10 & 5.8173e-12 & 1.5743e-13&4.8850e-15 & 1.1102e-15 
            \\
\hline
$O^l$ & -&  5.913 & 5.599 & 5.369 & 5.208 & 5.01 & 2.138
        \\
\hline
\end{tabular}
}
\caption{Grid refinement analysis for the accuracy of the quartic WENO spline and the classical quartic spline using the infinity norm. The original data has been sampled from the function in (\ref{pol}). We show the results for the functions $\psi_s, \psi_c, \psi_d$, presented in Subsection \ref{psi}, and the classical spline. The low-resolution nodes have been sampled with $m=2^l$ nodes and the high-resolution data with $11(m-1)+m$ points.}\label{tabla_quartic}
\end{center}
\end{table}

\begin{table}[!ht]
\begin{center}
%Programas en: /Users/j/Documents/articulos_a_medias/Levin/bsplines_WENO_pph/graficas_disc_articulo.m
\resizebox{15cm}{!} {
\begin{tabular}{|c|c|c|c|c|c|c|c|c|c|c|c|c|}
\hline $m=2^l$ & $2^4$ &  $2^5$           & $2^6$          &  $2^7$       & $2^8$ &  $2^9$   & $2^{10}$ %& $2^{11}$ & $2^{12}$
            \\
\hline
Error WENO-$\psi_s$ quintic spline ($E^l$)&      3.0918e-04 &  9.0968e-07 &  3.0077e-09 &  9.5293e-12 &  3.9302e-14 &  1.5543e-15 &  8.8818e-16
            \\
\hline
$O^l$ & -&       8.4089 &  8.2406 &  8.3021 &  7.9216 &  4.6603 &  8.0735e-01
        \\
\hline
Error WENO-$\psi_c$ quintic spline ($E^l$)&   2.0187e-05 & 2.0353e-08 & 1.4261e-10 & 1.9957e-12 & 3.1086e-14 & 1.5543e-15 & 1.1102e-15 
            \\
\hline
$O^l$ & -&   9.954 & 7.157 & 6.159 & 6.005 & 4.322 & 0.4854   
        \\      
\hline  
Error WENO-$\psi_d$ quintic spline ($E^l$)& 2.0204e-05 & 2.0353e-08 & 1.4261e-10 & 1.9957e-12 & 3.1086e-14 & 1.5543e-15 & 1.1102e-15   
            \\
\hline
$O^l$ & -&  9.955 & 7.157 & 6.159 & 6.005 & 4.322 & 0.4854 
        \\
\hline
\hline
Error class. quintic spline ($E^l$)&       7.2832e-07 &  9.3475e-09 &  1.3269e-10 &  1.9780e-12 &  3.1308e-14 &  1.3323e-15 &  1.1102e-15
            \\
\hline
$O^l$ & -&       6.2838 &  6.1385 &  6.0679 &  5.9813 &  4.5546 &  2.6303e-01
        \\
\hline
\end{tabular}
}
\caption{Grid refinement analysis for the accuracy of the quintic WENO spline and the classical quintic spline using the infinity norm. The original data has been sampled from the function in (\ref{pol}). We show the results for the functions $\psi_s, \psi_c, \psi_d$, presented in Subsection \ref{psi}, and the classical spline. The low-resolution nodes have been sampled with $m=2^l$ nodes and the high-resolution data with $10(m-1)+m$ points.}\label{tabla_quintic}
\end{center}
\end{table}

\subsection{Accuracy close to the discontinuity}

It is also interesting to check the numerical accuracy of the new method close to the discontinuity, as stated in Preposition \ref{ordenh}. We consider the function in (\ref{sin})
\begin{equation}\label{sin}
f(x)=\left\{\begin{array}{ll}
\cos(x-0.5), \quad 0\le x\le 0.5,\\
\sin(x), \quad 0.5<x\le 1.
\end{array}\right.
\end{equation}
We check the accuracy in $(0.5, 1]$ but avoiding the interval that contains the discontinuity close to $x=0.5$. The reason is that we know that the approximation in the interval that contains the discontinuity presents the order of accuracy is $O(h^0)$. To do so, as before, we use the infinity norm of the error in (\ref{linf}) and the formula for the numerical order in (\ref{orden}). The results for the quadratic spline are presented in table \ref{quadratic_cerca_disc}, for the cubic in table \ref{cubic_cerca_disc}, for the quartic in table \ref{quartic_cerca_disc} and for the quintic in table \ref{quintic_cerca_disc}. As before, all the tables show the results for the functions $\psi_s, \psi_c, \psi_d$. All of them show order $O(h)$ in the infinity norm, as proved in Preposition \ref{ordenh}. This order of approximation corresponds to the intervals that are adjacent to the one that contains the discontinuity.

%\begin{landscape}
\begin{table}[!ht]
\begin{center}
%Programas en: /Users/j/Documents/articulos_a_medias/Levin/bsplines_WENO_pph/graficas_disc_articulo.m
\resizebox{18cm}{!} {
\begin{tabular}{|c|c|c|c|c|c|c|c|c|c|c|c|c|c|c|c|c|c|c|c|c|c|c|c|c|c|c|}
%\hline $i$ & $4$  & $5$      & $6$     & $7$    & $8$ & $9$    & $10$
%\\
%\hline $n_i$ & 16 & 32& 64      & 128     & 256    & 512 & 1024
\hline $m=2^l$ & $2^4$ &  $2^5$           & $2^6$          &  $2^7$       & $2^8$ &  $2^9$   & $2^{10}$ & $2^{11}$ & $2^{12}$ & $2^{13}$
      \\
\hline
Error WENO-$\psi_s$ quadratic spline ($E^l$)&   5.4615e-02 & 2.7978e-02 & 1.3927e-02 & 6.9191e-03 & 3.4449e-03 & 1.7184e-03 & 8.5812e-04 & 4.2879e-04 & 2.1432e-04 & 1.0714e-04
            \\
\hline
$O^l$ & - & 0.965 & 1.006 & 1.009 & 1.006 & 1.003& 1.002 & 1.001 &    1 &    1
        \\
\hline        
Error WENO-$\psi_c$ quadratic spline ($E^l$)&  4.3329e-02 & 2.9412e-02 & 1.8279e-02 & 1.0507e-02 & 5.6976e-03 & 2.9776e-03 & 1.5237e-03 & 7.7091e-04 & 3.8777e-04 & 1.9447e-04
            \\
\hline
$O^l$ &- & 0.5589 & 0.6862 & 0.7989 & 0.8829 & 0.9362 & 0.9666 & 0.9829 & 0.9914 & 0.9956
        \\        
\hline
Error WENO-$\psi_d$ quadratic spline ($E^l$)&   6.0303e-02 & 2.7966e-02 & 1.3823e-02 & 6.8839e-03 & 3.4350e-03 & 1.7158e-03 & 8.5745e-04 & 4.2862e-04 & 2.1428e-04 & 1.0713e-04
            \\
\hline
$O^l$ &- & 1.109 & 1.017 & 1.006 & 1.003 & 1.001 & 1.001 &    1 &    1 &    1
        \\
\hline
\end{tabular}
}
\caption{Grid refinement analysis in the infinity norm close to the discontinuity of the function shown in (\ref{sin}). In this case, we present the results of the quadratic WENO algorithm for the functions $\psi_s, \psi_c, \psi_d$, presented in Subsection \ref{psi}.
}\label{quadratic_cerca_disc}
\end{center}
\end{table}

\begin{table}[!ht]
\begin{center}
%Programas en: /Users/j/Documents/articulos_a_medias/Levin/bsplines_WENO_pph/graficas_disc_articulo.m
\resizebox{15cm}{!} {
\begin{tabular}{|c|c|c|c|c|c|c|c|c|c|c|c|c|}
\hline $m=2^l$ & $2^4$ &  $2^5$           & $2^6$          &  $2^7$       & $2^8$ &  $2^9$   & $2^{10}$ & $2^{11}$ & $2^{12}$ & $2^{13}$
            \\
\hline
Error WENO-$\psi_s$ cubic spline ($E^l$)&  6.4613e-02 & 3.2117e-02 & 1.5988e-02 & 7.9734e-03 & 3.9812e-03 & 1.9892e-03 & 9.9425e-04 & 4.9704e-04 & 2.4850e-04 & 1.2424e-04
            \\
\hline
$O^l$ &- & 1.008 & 1.006 & 1.004 & 1.002 & 1.001 & 1.001 &    1 &    1 &    1
        \\
\hline
Error WENO-$\psi_c$ cubic spline ($E^l$)&  5.7789e-02 & 2.8948e-02 & 1.4497e-02 & 7.2559e-03 & 3.6302e-03 & 1.8157e-03 & 9.0799e-04 & 4.5404e-04 & 2.2703e-04 & 1.1352e-04
            \\
\hline
$O^l$ &- & 0.9973 & 0.9978 & 0.9985 & 0.9991 & 0.9995 & 0.9998 & 0.9999 & 0.9999 &    1
        \\                
\hline
Error WENO-$\psi_d$ cubic spline ($E^l$)&   6.4763e-02 & 3.2233e-02 & 1.6013e-02 & 7.9793e-03 & 3.9826e-03 & 1.9896e-03 & 9.9434e-04 & 4.9706e-04 & 2.4850e-04 & 1.2424e-04
            \\
\hline
$O^l$ & - & 1.007 & 1.009 & 1.005 & 1.003 & 1.001& 1.001 &    1 &    1 &    1
        \\
\hline
\end{tabular}
}
\caption{Grid refinement analysis in the infinity norm close to the discontinuity of the function shown in (\ref{sin}). In this case, we present the results of the cubic WENO algorithm for the functions $\psi_s, \psi_c, \psi_d$, presented in Subsection \ref{psi}.}\label{cubic_cerca_disc}
\end{center}
\end{table}

\begin{table}[!ht]
\begin{center}
%Programas en: /Users/j/Documents/articulos_a_medias/Levin/bsplines_WENO_pph/graficas_disc_articulo.m
\resizebox{15cm}{!} {
\begin{tabular}{|c|c|c|c|c|c|c|c|c|c|c|c|c|}
\hline $m=2^l$ & $2^4$ &  $2^5$           & $2^6$          &  $2^7$       & $2^8$ &  $2^9$   & $2^{10}$ & $2^{11}$ & $2^{12}$& $2^{13}$
            \\
\hline
Error WENO-$\psi_s$ quartic spline ($E^l$)& 7.0655e-02 & 5.1394e-02 & 3.9749e-02 & 3.1339e-02 & 2.3823e-02 & 1.6675e-02 & 1.0554e-02 & 6.1095e-03 & 3.3201e-03 & 1.7359e-03
            \\
\hline
$O^l$ &- & 0.4592 & 0.3707 & 0.3429 & 0.3956 & 0.5147& 0.6599 & 0.7886 & 0.8799 & 0.9355
        \\
\hline
Error WENO-$\psi_c$ quartic spline ($E^l$)&  7.0655e-02 & 5.1394e-02 & 3.9749e-02 & 3.1339e-02 & 2.3823e-02 & 1.6675e-02 & 1.0554e-02 & 6.1095e-03 & 3.3201e-03 & 1.7359e-03
            \\
\hline
$O^l$ & - & 0.4592 & 0.3707 & 0.3429 & 0.3956 & 0.5147& 0.6599 & 0.7886 & 0.8799 & 0.9355
        \\        
\hline
Error WENO-$\psi_d$ quartic spline ($E^l$)&  7.0655e-02 & 5.1394e-02 & 3.9749e-02 & 3.1339e-02 & 2.3823e-02 & 1.6675e-02 & 1.0554e-02 & 6.1095e-03 & 3.3201e-03 & 1.7359e-03
            \\
\hline
$O^l$ &- & 0.4592 & 0.3707 & 0.3429 & 0.3956 & 0.5147& 0.6599 & 0.7886 & 0.8799 & 0.9355
        \\
\hline
\end{tabular}
}
\caption{Grid refinement analysis in the infinity norm close to the discontinuity of the function shown in (\ref{sin}). In this case, we present the results of the quartic WENO algorithm for the functions $\psi_s, \psi_c, \psi_d$, presented in Subsection \ref{psi}.}\label{quartic_cerca_disc}
\end{center}
\end{table}

\begin{table}[!ht]
\begin{center}
%Programas en:/Users/j/Documents/articulos_a_medias/Levin/bsplines_WENO_pph/graficas_disc_articulo.m
\resizebox{15cm}{!} {
\begin{tabular}{|c|c|c|c|c|c|c|c|c|c|c|c|c|}
\hline $m=2^l$ & $2^4$ &  $2^5$           & $2^6$          &  $2^7$       & $2^8$ &  $2^9$   & $2^{10}$ & $2^{11}$ & $2^{12}$& $2^{13}$
            \\
\hline
Error WENO-$\psi_s$ quintic spline ($E^l$)& 1.0267e-01 & 5.6743e-02 & 2.8517e-02 & 1.4064e-02 & 6.9536e-03 & 3.4536e-03 & 1.7206e-03 & 8.5866e-04 & 4.2892e-04 & 2.1436e-04
            \\
\hline
$O^l$ & - & 0.85548 & 0.99262 & 1.0198 & 1.0162 & 1.0097& 1.0052 & 1.0027 & 1.0014 & 1.0007
        \\
\hline
Error WENO-$\psi_c$ quintic spline ($E^l$)&  8.8985e-02 & 6.5004e-02 & 4.7851e-02 & 3.3716e-02 & 2.1894e-02 & 1.3017e-02 & 7.2131e-03 & 3.8165e-03 & 1.9660e-03 & 9.9813e-04
            \\
\hline
$O^l$ &- & 0.45304 & 0.44199 & 0.50511 & 0.62287 & 0.75013& 0.85174 & 0.91836 & 0.95702 & 0.97793
        \\      
\hline  
Error WENO-$\psi_d$ quintic spline ($E^l$)& 1.1271e-01 & 5.5344e-02 & 2.7521e-02 & 1.3738e-02 & 6.8626e-03 & 3.4297e-03 & 1.7144e-03 & 8.5712e-04 & 4.2853e-04 & 2.1426e-04
            \\
\hline
$O^l$ &- & 1.0261 & 1.0079 & 1.0024 & 1.0013 & 1.0007& 1.0003 & 1.0002 & 1.0001 &    1
        \\
\hline
\end{tabular}
}
\caption{Grid refinement analysis in the infinity norm close to the discontinuity of the function shown in (\ref{sin}). In this case, we present the results of the quintic WENO algorithm for the functions $\psi_s, \psi_c, \psi_d$, presented in Subsection \ref{psi}.}\label{quintic_cerca_disc}
\end{center}
\end{table}

%%\begin{landscape}
%\begin{table}[!ht]
%\begin{center}
%%Programas en: /Users/j/Documents/articulos_a_medias/Levin/bsplines_WENO_pph/refinamiento_malla_articulo.m
%%Programas en: /Users/j/Documents/articulos_a_medias/Levin/bsplines_WENO_pph/weno_spline.m
%\resizebox{15cm}{!} {
%\begin{tabular}{|c|c|c|c|c|c|c|c|c|c|c|c|c|c|c|c|c|c|c|c|c|c|c|c|c|c|c|}
%\hline $i$ & $1$  & $2$      & $3$     & $4$    & $5$
%\\
%\hline $n_i$ & $101$  & $201$      & $401$     & $801$    & $1601$% & $2^9$  & $2^{10}$ &$2^{11}$ &$2^{12}$&$2^{13}$%&$2^{14}$&$2^{15}$
%      \\
%\hline
%Error ($E_i^{\infty}$)&   0.0088  &  0.0044 &   0.0022 &   0.0011 &  5.4843e-04
%      \\
%\hline
%$O_i$&- & 1.0013  & 1.0013  & 1.0013  & 1.0054
%    \\
%\hline
%\end{tabular}
%}
%\caption{Grid refinement analysis in the infinity norm for the function shown in (\ref{sin}).
%}\label{tabla}
%\end{center}
%\end{table}

\subsection{Adaption to the presence of discontinuities: Gibbs phenomenon}

In this subsection, we analyze the behavior of the algorithms  in terms of their adaption close the discontinuity. Figure \ref{figuras_1D} shows the results of approximating the function in (\ref{sin}) close to the discontinuity. To obtain these results, we have started from a sampling of the function in (\ref{sin}) of 400 initial points. Then we obtained approximations at the position of these initial points, plus 11 intermediate points for the quadratic and quartic spline, or 10 intermediate points for the cubic and quintic spline. As before, we present the results for the quadratic, cubic, quartic, and quintic WENO B-spline-based algorithm for the functions $\psi_s, \psi_c, \psi_d$, presented in Subsection \ref{psi}, and the respective classical splines. It is easy to appreciate how the three versions of the WENO B-spline-based algorithm presented managed to reduce the oscillation close to the discontinuity that the classical spline presents.  We can also see that each function $\psi_s, \psi_c$ or $\psi_d$ allows the WENO algorithm to present different degrees of sharpness close to the discontinuity, $\psi_d$ being the sharpest one. For the cubic spline, we can see that the three functions $\psi_s, \psi_c$, and $\psi_d$ behave very similarly close to the discontinuity. For the quintic spline $\psi_s$ and $\psi_c$ also behave very similarly. It is important to mention that at the interval that contains the discontinuity all the algorithms present smearing, so it is not possible to control the error.
\begin{figure}[!ht]
%/Users/j/Documents/articulos_a_medias/Levin/bsplines_WENO_pph/graficas_disc_articulo.m
\centerline{\psfig{figure=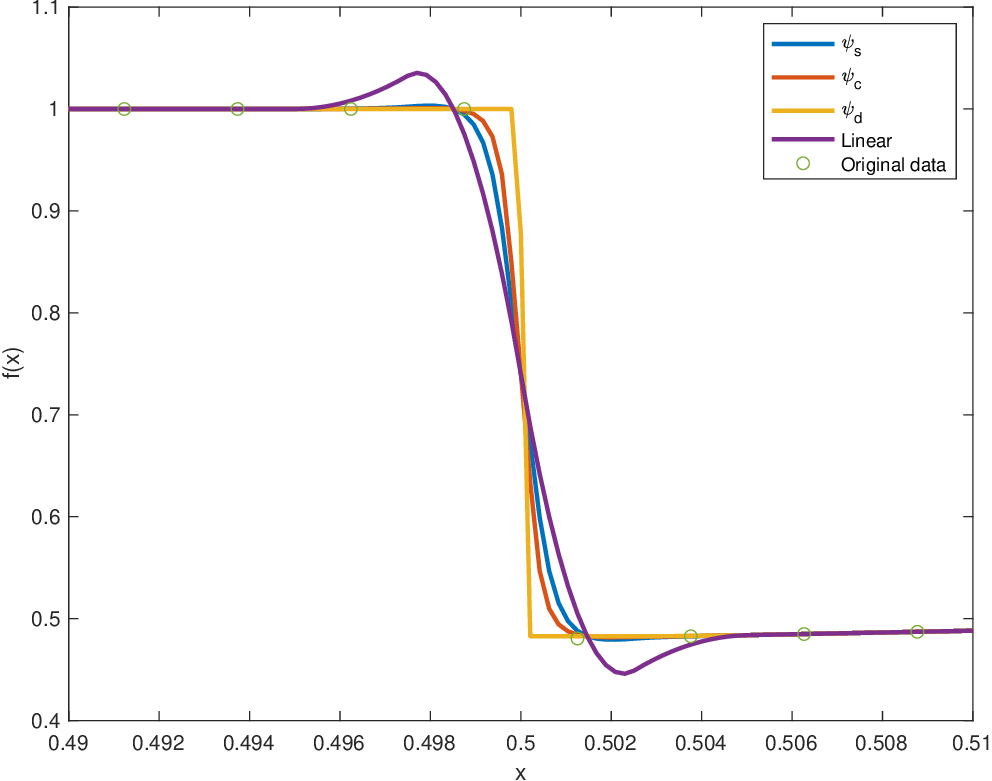,height=4cm}
\psfig{figure=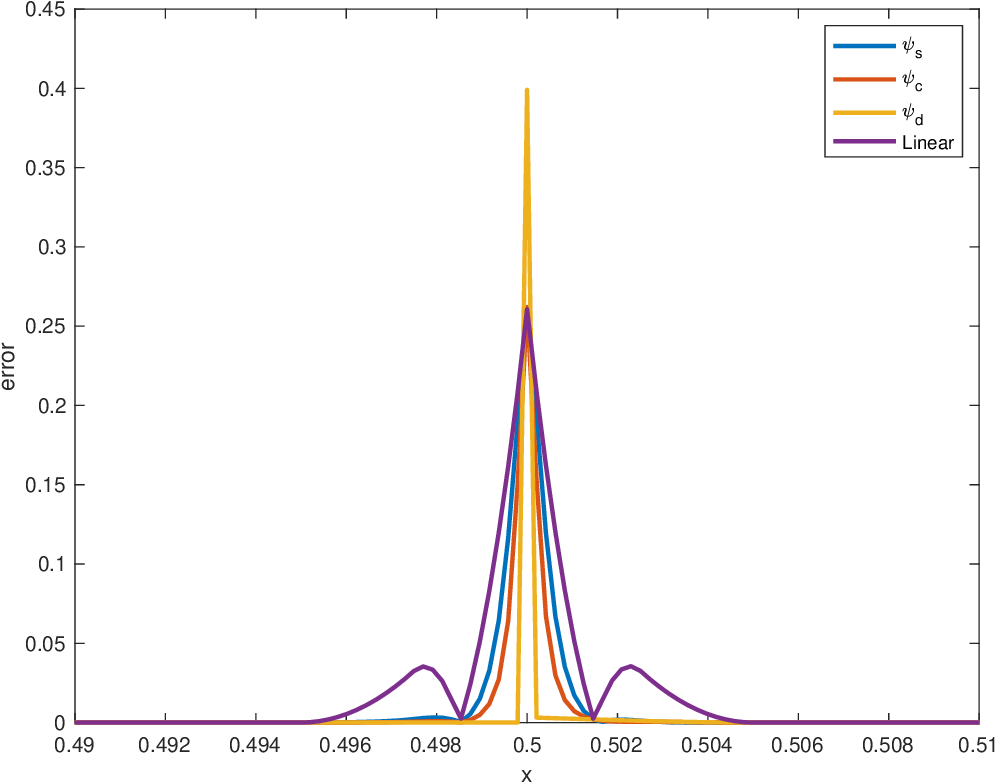,height=4cm}}
\centerline{\psfig{figure=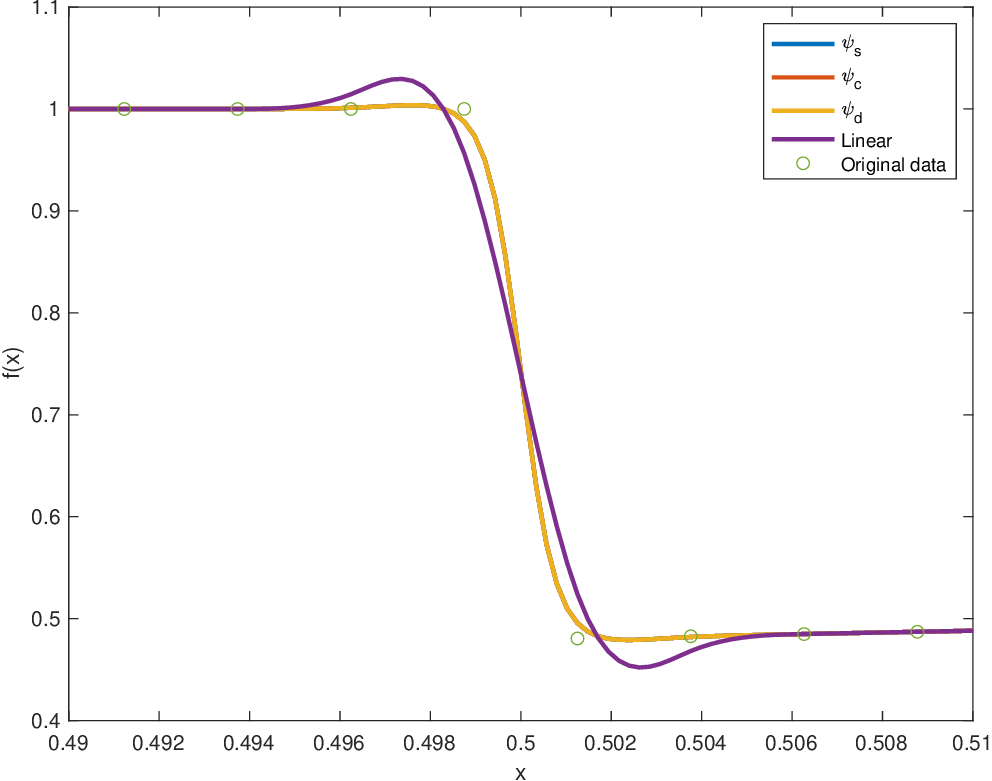,height=4cm}
\psfig{figure=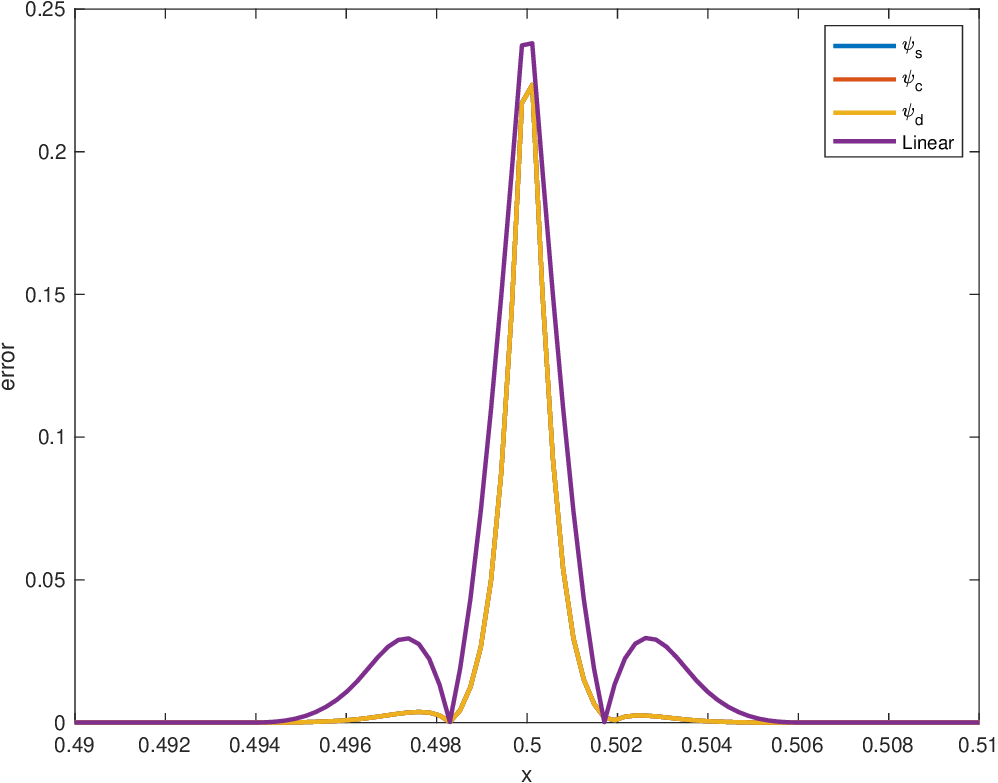,height=4cm}}
\centerline{\psfig{figure=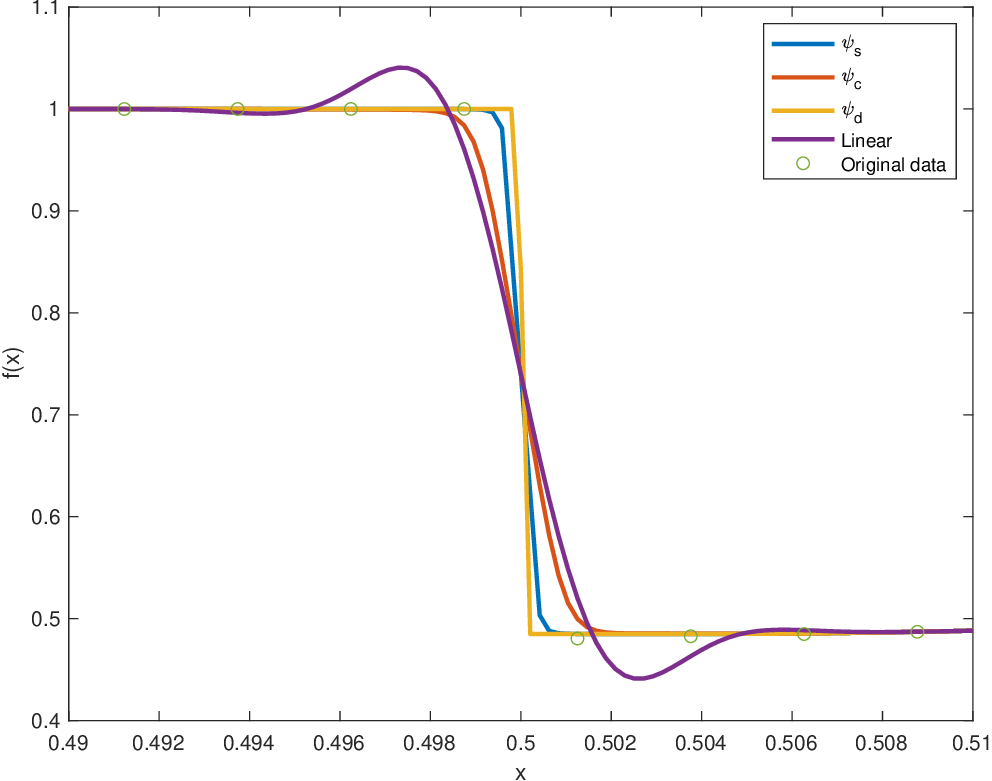,height=4cm}
\psfig{figure=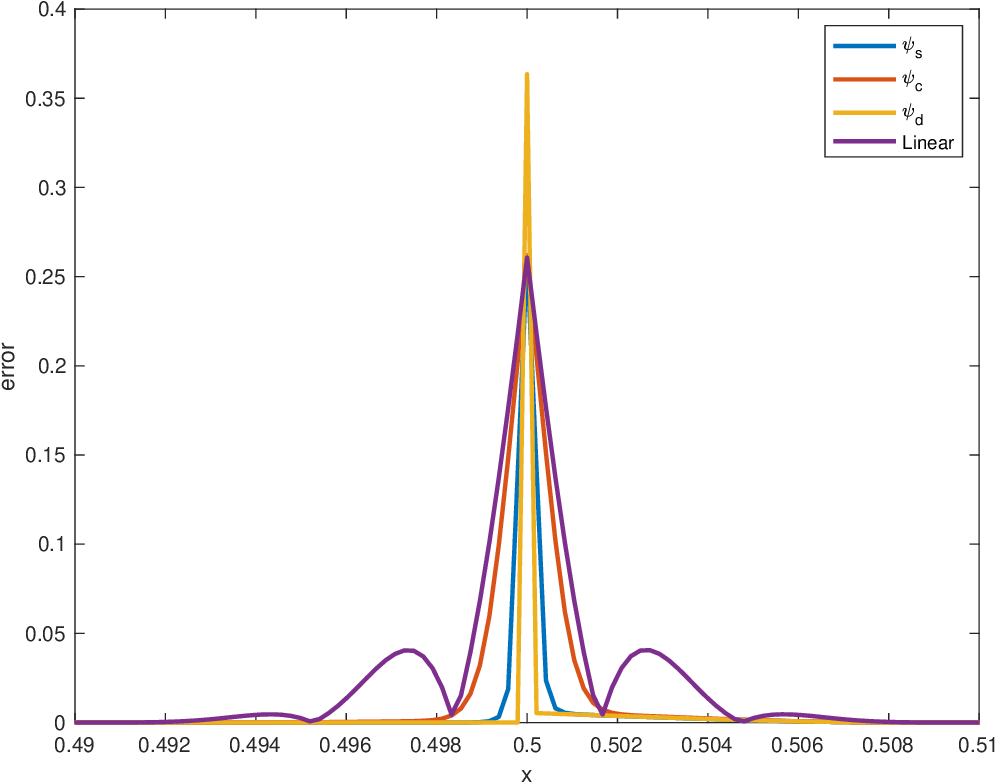,height=4cm}}
\centerline{\psfig{figure=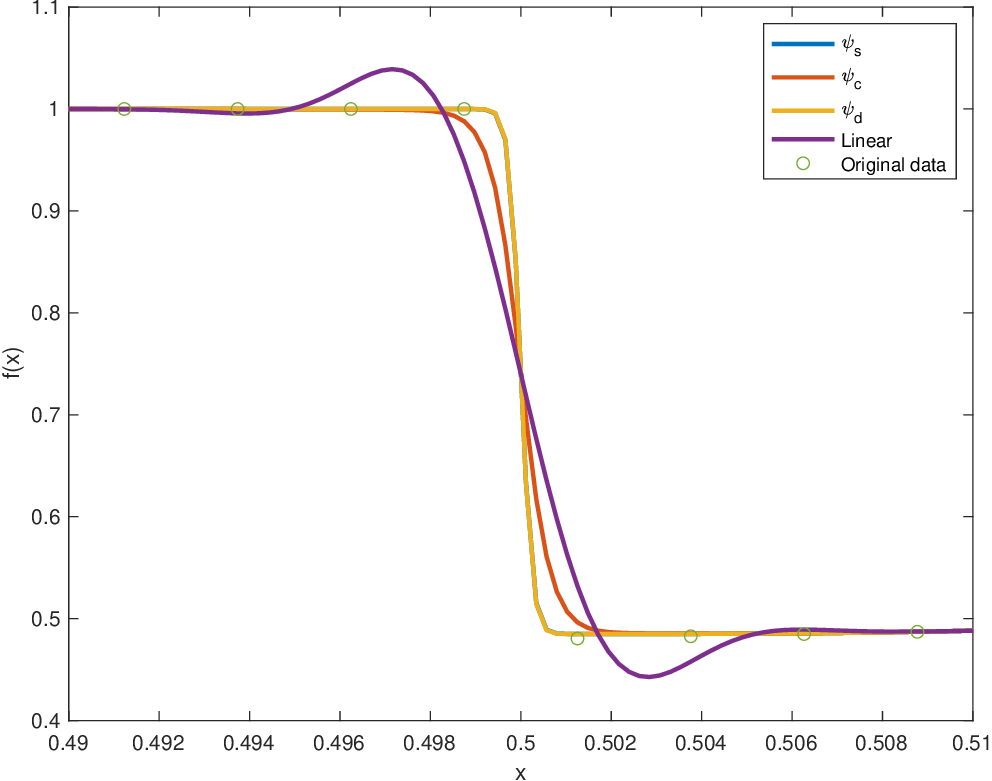,height=4cm}
\psfig{figure=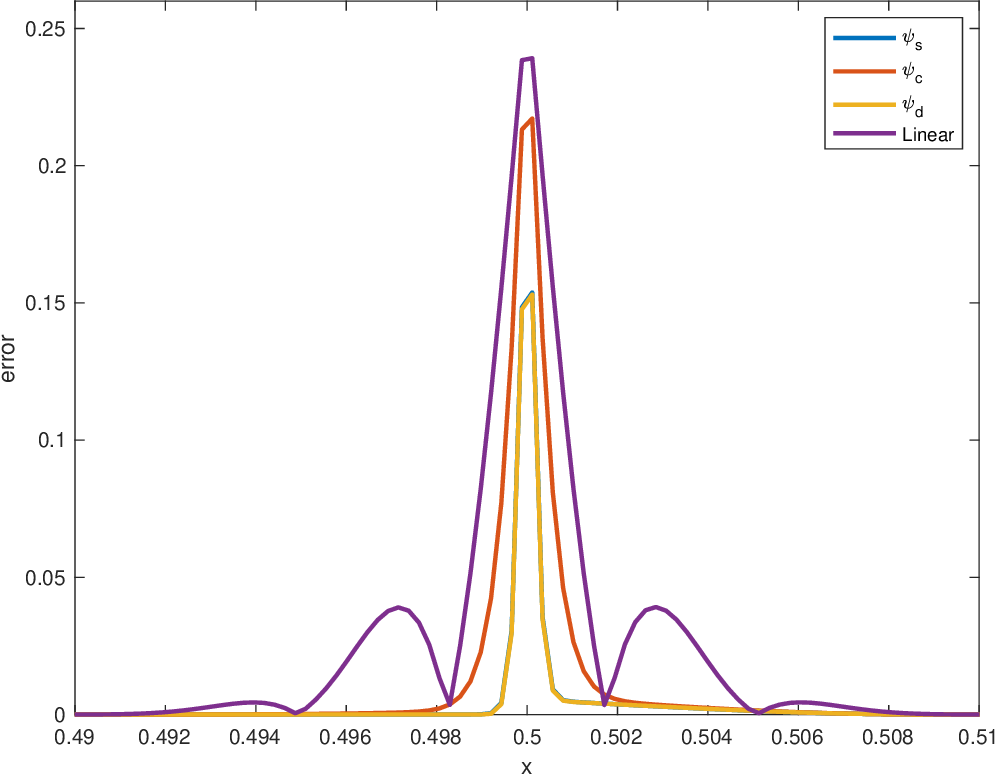,height=4cm}}
\caption{The column to the left presents the limit function obtained for $\psi_s,\psi_c, \psi_d$ and the linear algorithm. The column to the right presents the error obtained by each of them. The first row presents the result of the quadratic algorithm, the second presents the results of the cubic, and so on.}\label{figuras_1D}
\end{figure}

\subsection{Experiments for cubic splines and bi-variate data using tensor product}
In this section, we present the results that we obtain when extending the one-dimensional algorithms presented in previous sections to several dimensions. The process is described in Section \ref{dimensions}. Figure \ref{figuras_2D} shows the results obtained when approximating the function

\begin{equation}\label{2D}
f(x,y)=\left\{\begin{array}{ll}
\cos(xy), \quad (x-0.5)^2+(y-0.5)^2\le r^2,\\
\sin(xy), \quad (x-0.5)^2+(y-0.5)^2>r^2,
\end{array}\right.
\end{equation}
with $r=\frac{1}{4}$. The results are shown in Figure \ref{2D}. We can see how the nonlinear algorithms manage to reduce the oscillations of the classical approach close to the discontinuity.

\begin{figure}[!ht]
%/Users/j/Documents/articulos_a_medias/Levin/bsplines_WENO_pph/graficas_disc_articulo.m
\centerline{\psfig{figure=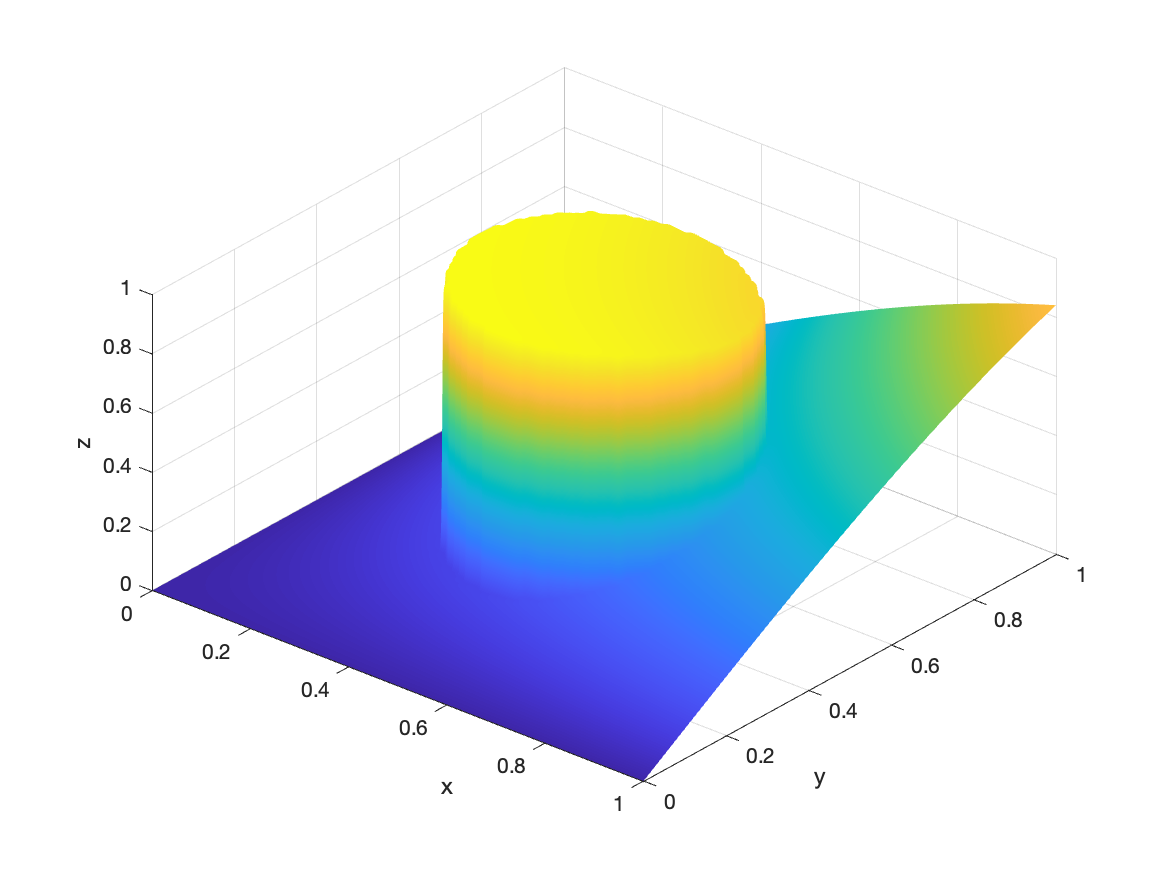,height=4cm}
\psfig{figure=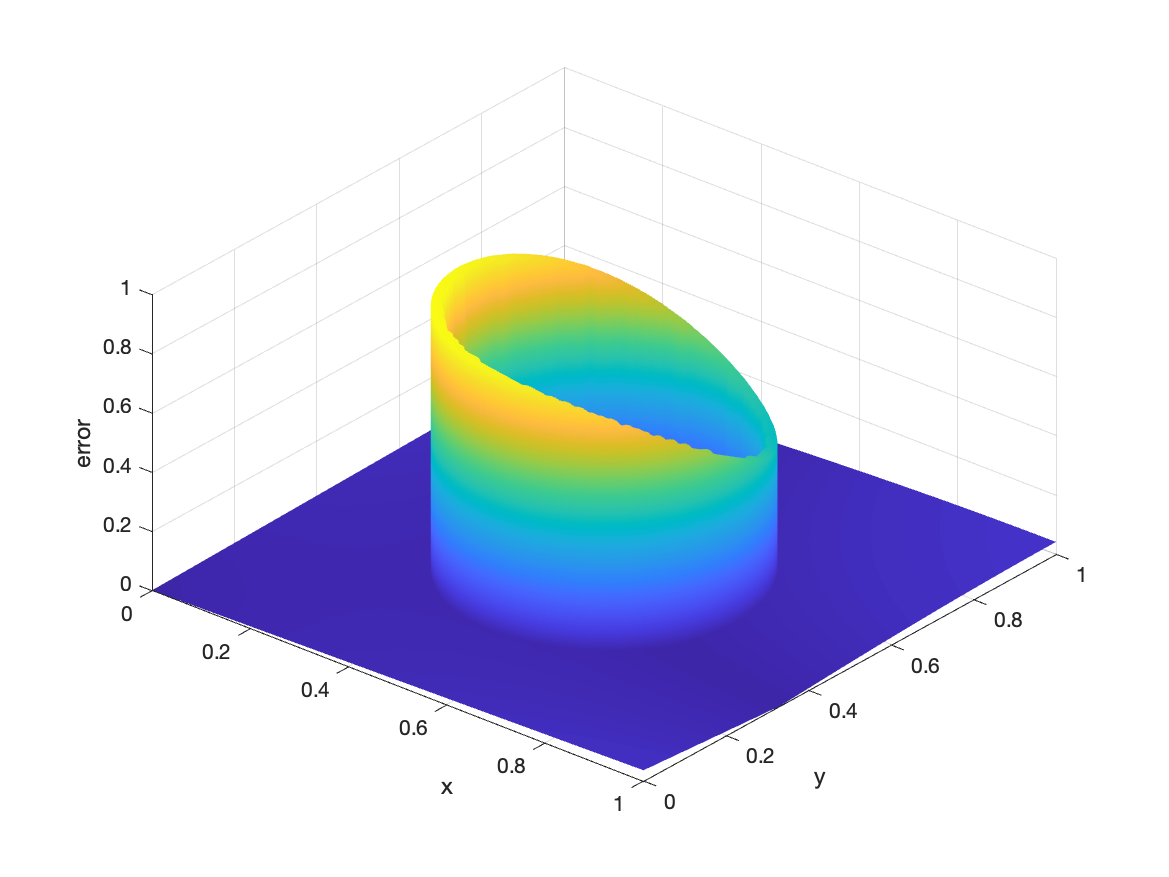,height=4cm}}
\centerline{\psfig{figure=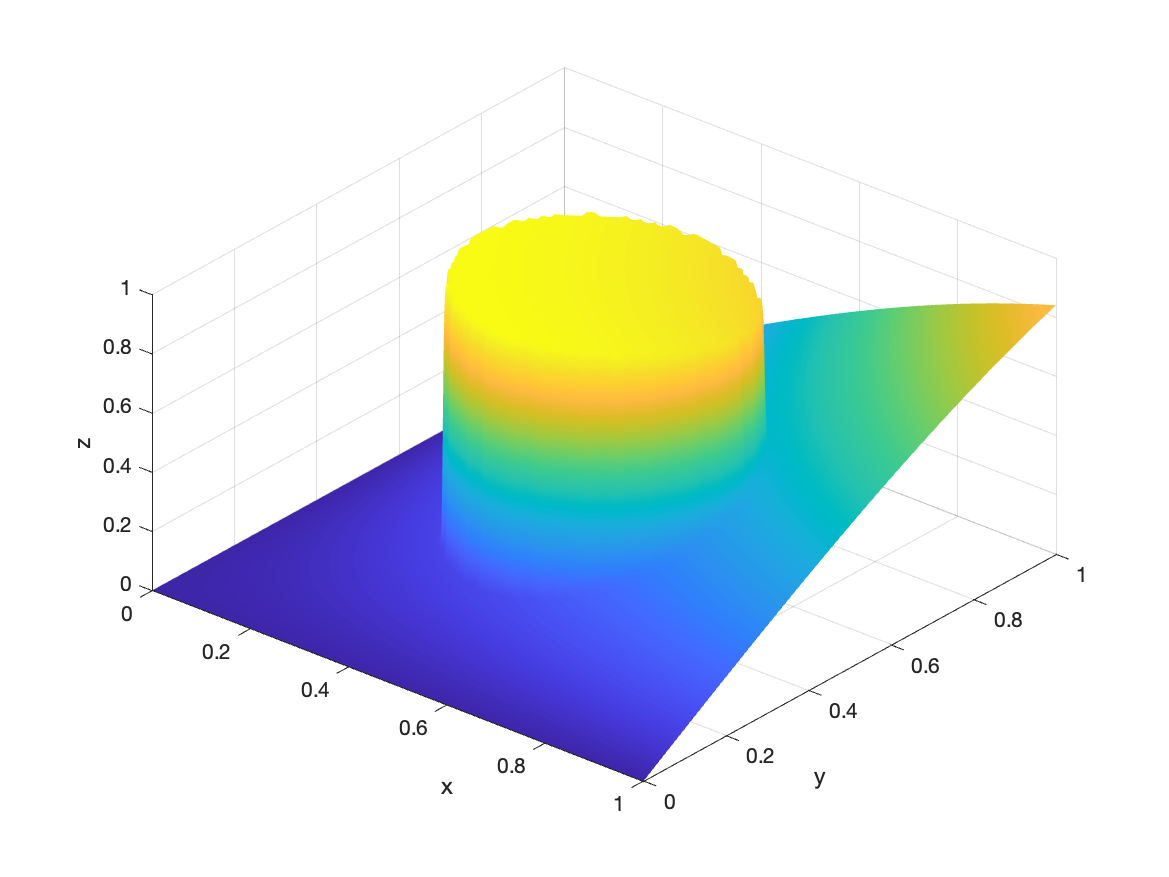,height=4cm}
\psfig{figure=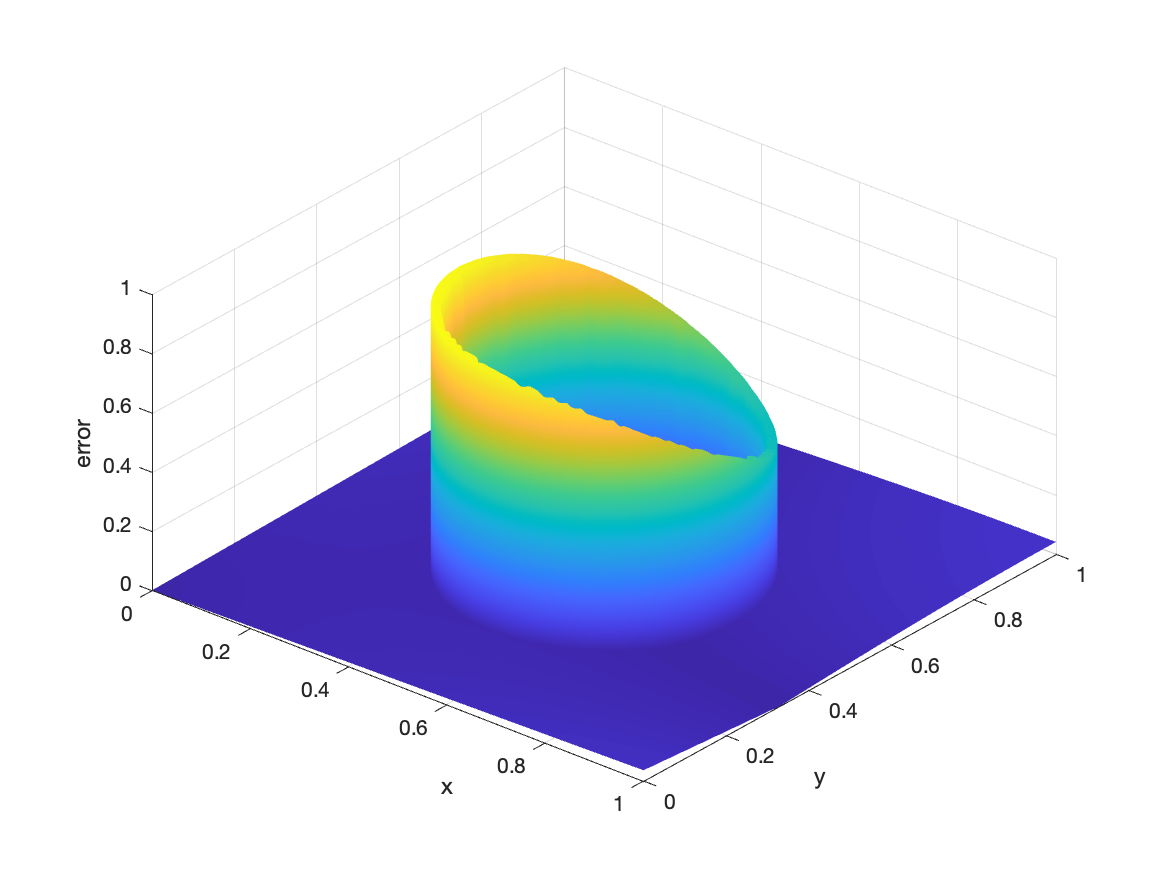,height=4cm}}
\centerline{\psfig{figure=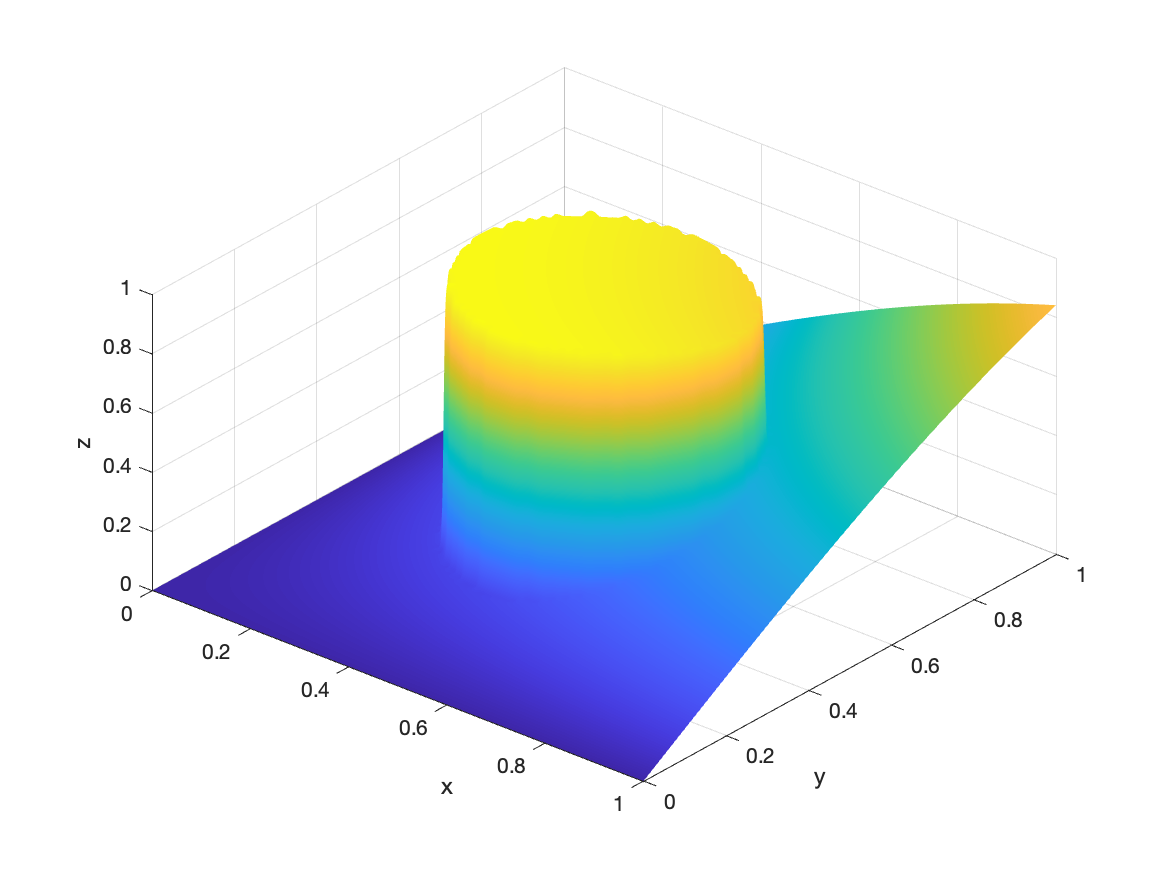,height=4cm}
\psfig{figure=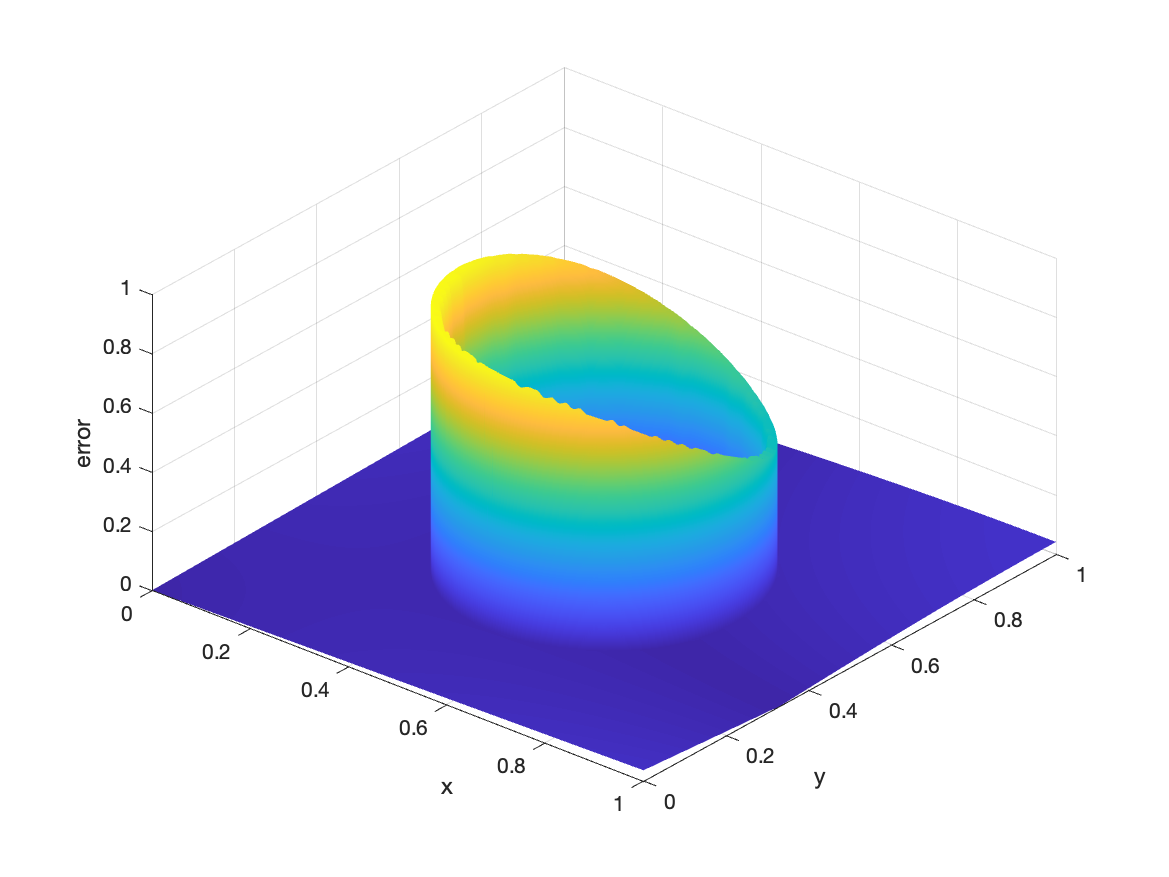,height=4cm}}
\centerline{\psfig{figure=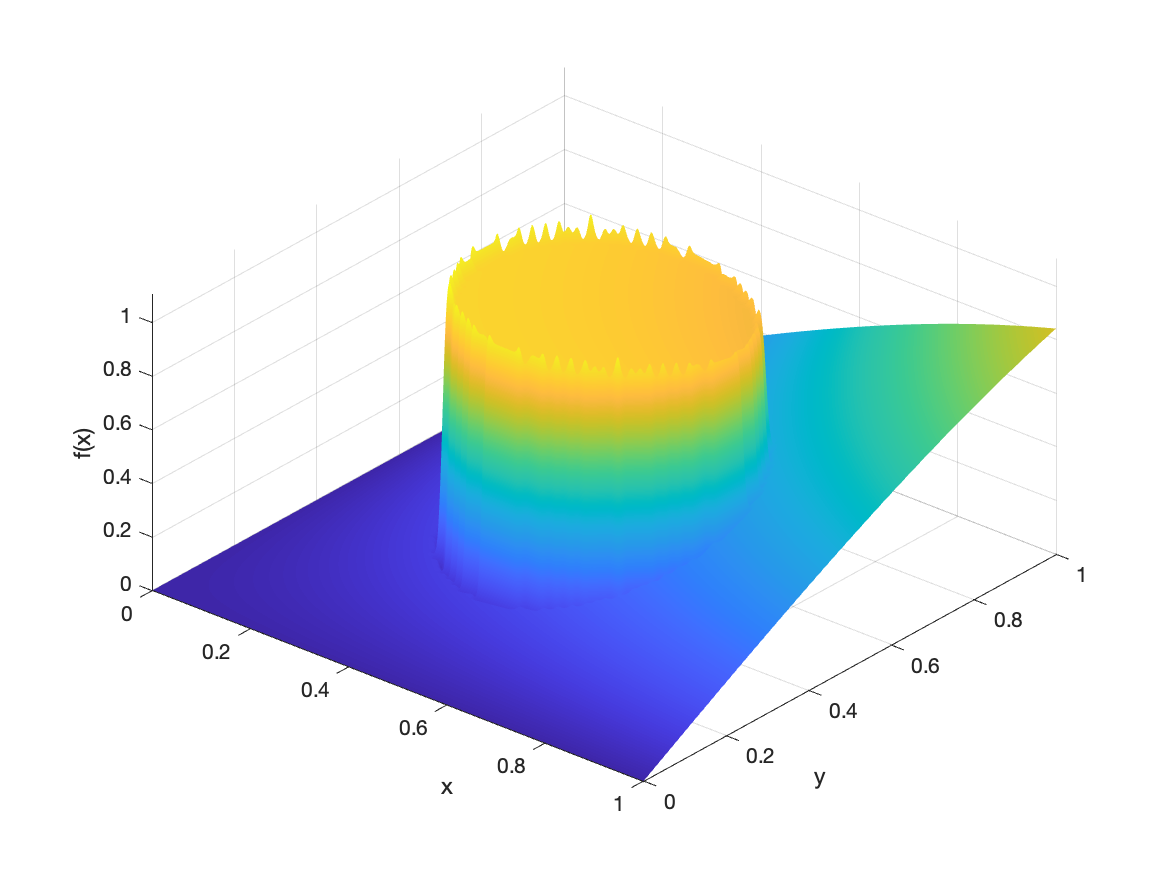,height=4cm}
\psfig{figure=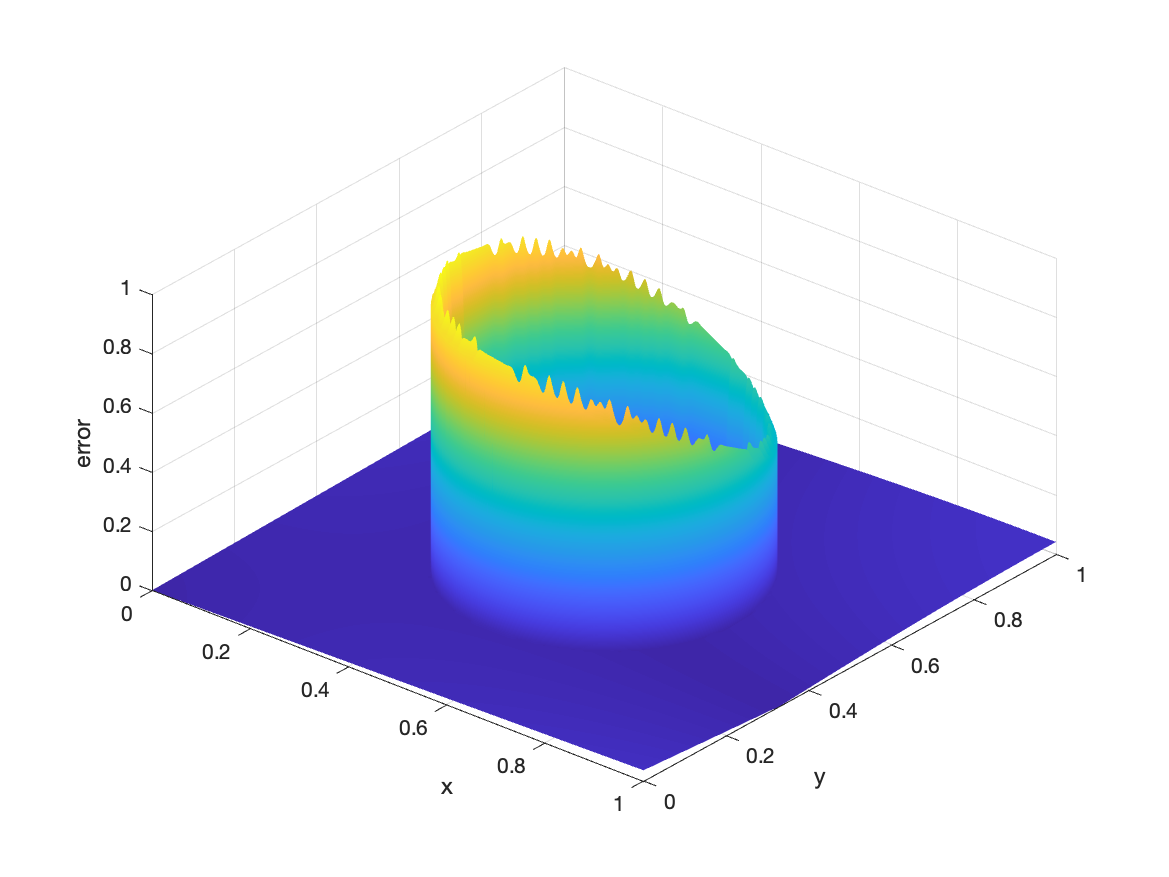,height=4cm}}
\caption{Left column, approximation of the function in (\ref{2D}) obtained for $\psi_s$ (first row), $\psi_c$ (second row), $\psi_d$ (third row) and the linear algorithm (fourth row). Right column, the distribution of the error obtained by each of them.}\label{figuras_2D}
\end{figure}

\subsection{Experiments for cubic splines and three-variate data using tensor product}

In this section, we present similar results to the ones obtained in the previous subsection, but applied to trivariate data. We have used again the tensor product approach described in Section \ref{dimensions}. Figure \ref{figuras_3D} shows the results obtained when approximating the function

\begin{equation}\label{3D}
f(x,y,z)=\left\{\begin{array}{ll}
e^{x+y+z}, \quad (x-0.5)^2+(y-0.5)^2+(z-0.5)^2\le r^2,\\
\cos(x+y+z), \quad (x-0.5)^2+(y-0.5)^2+(z-0.5)^2>r^2,
\end{array}\right.
\end{equation}
with $r=0.4$. In this case, we have chosen the function for the weights $\psi_d$, presented in Subsection \ref{psi}. To obtain the results shown in Figure \ref{3D}, the function in (\ref{3D}) has been sampled with $200\times200\times200$ points. The reconstruction has been obtained at the position of the initial values plus at two equidistant intermediate values between the original ones in each direction. In this case, we present the approximation of the trivariate function in (\ref{3D}) using a contour slice routine: we have obtained slices of the volumetric data in the $z$ direction and we have represented contour plots of the bivariate data obtained on each slice. We can see how the nonlinear algorithm manages to reduce the oscillations of the classical approach close to the discontinuity, that in this case is a surface, obtaining a sharp reconstruction close to it. The oscillations in the 3D reconstruction of the linear algorithm can be inferred from the high-frequency oscillations in the contour plots close to the discontinuity surface and from the presence of many small closed contours also close to the discontinuity surface.

\begin{figure}[!ht]
%/Users/j/Documents/articulos_a_medias/Levin/bsplines_WENO_pph/graficas_disc_articulo.m
\centerline{\psfig{figure=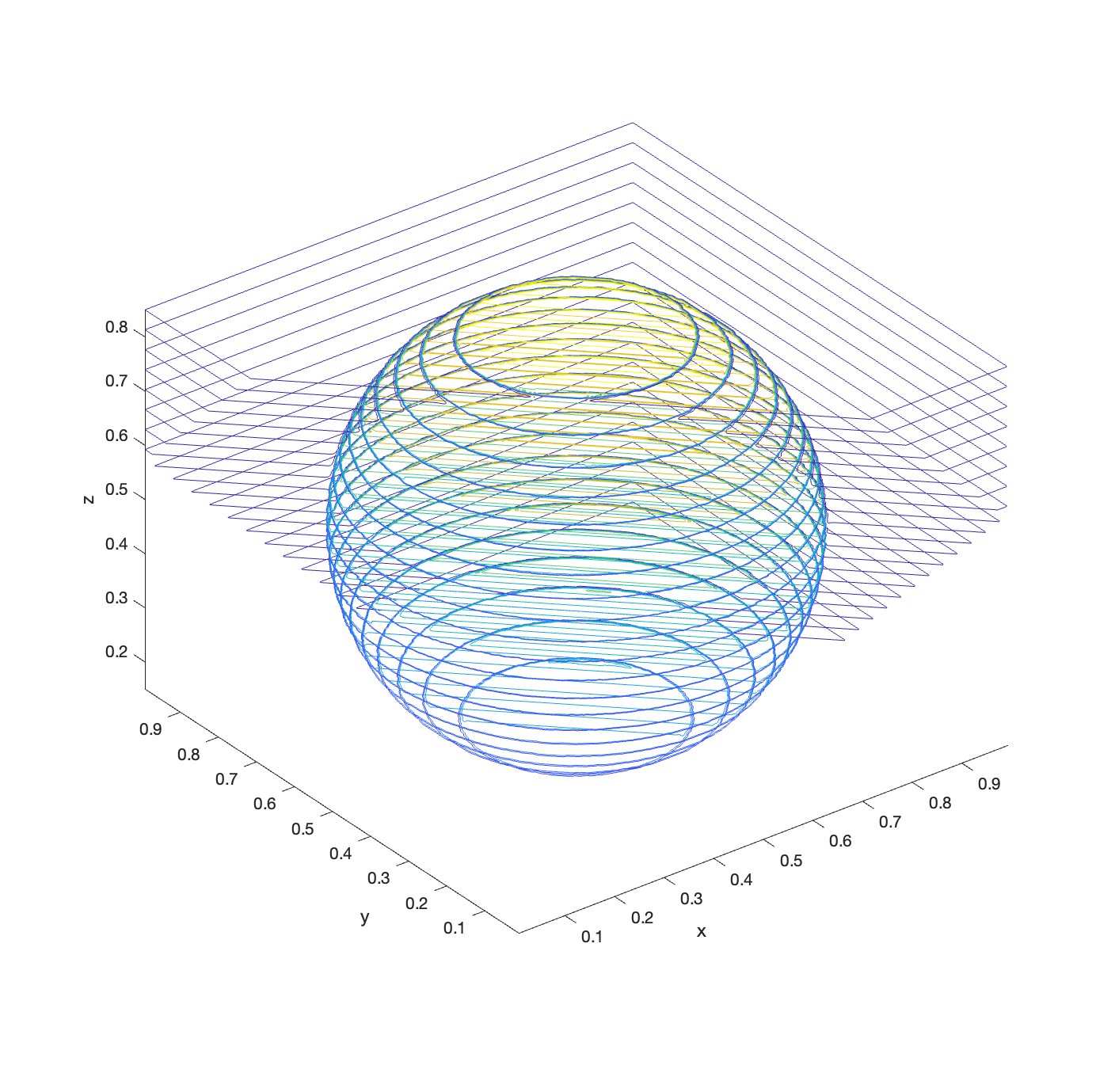,height=8cm}
\psfig{figure=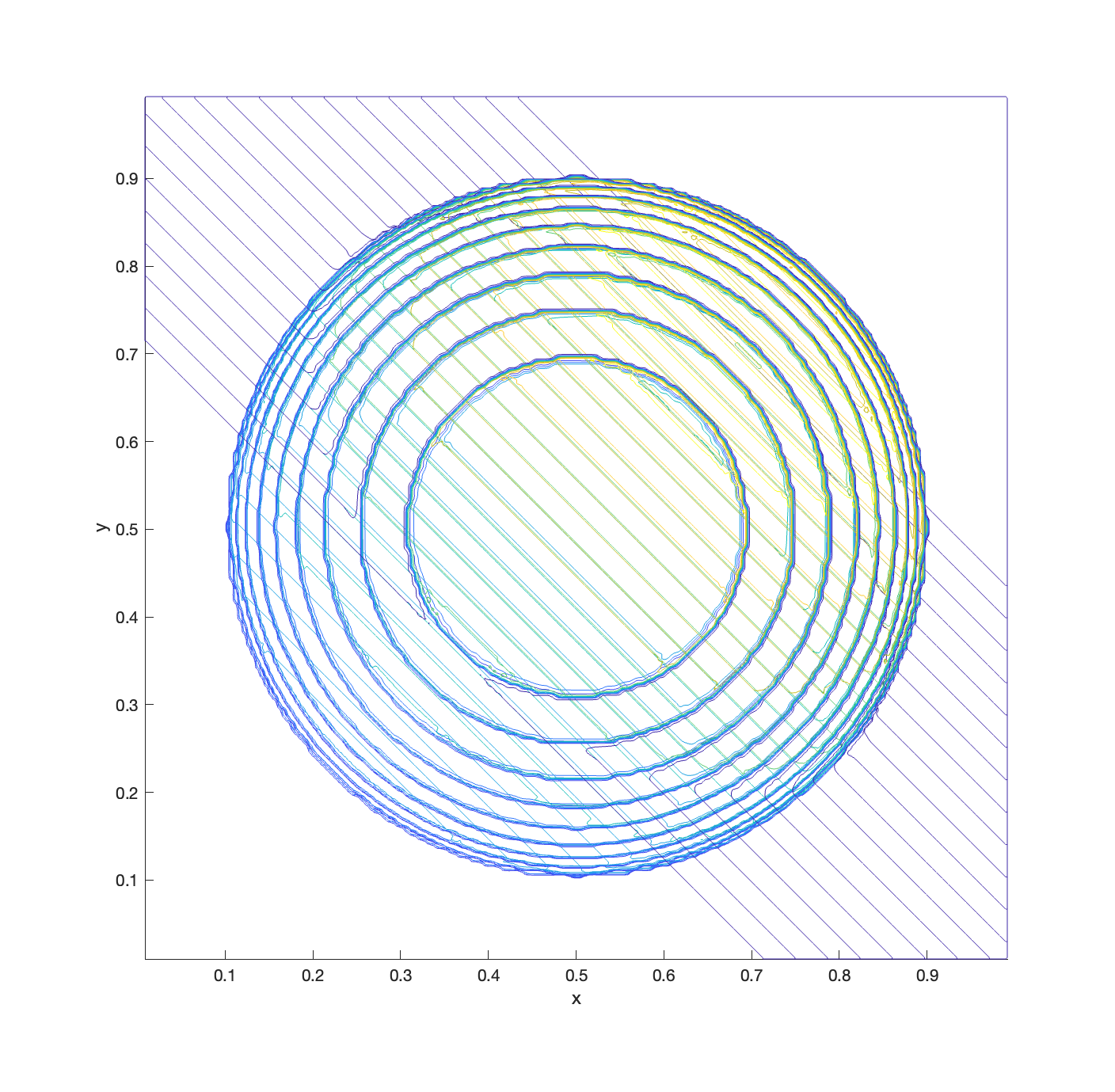,height=8cm}}
\centerline{\psfig{figure=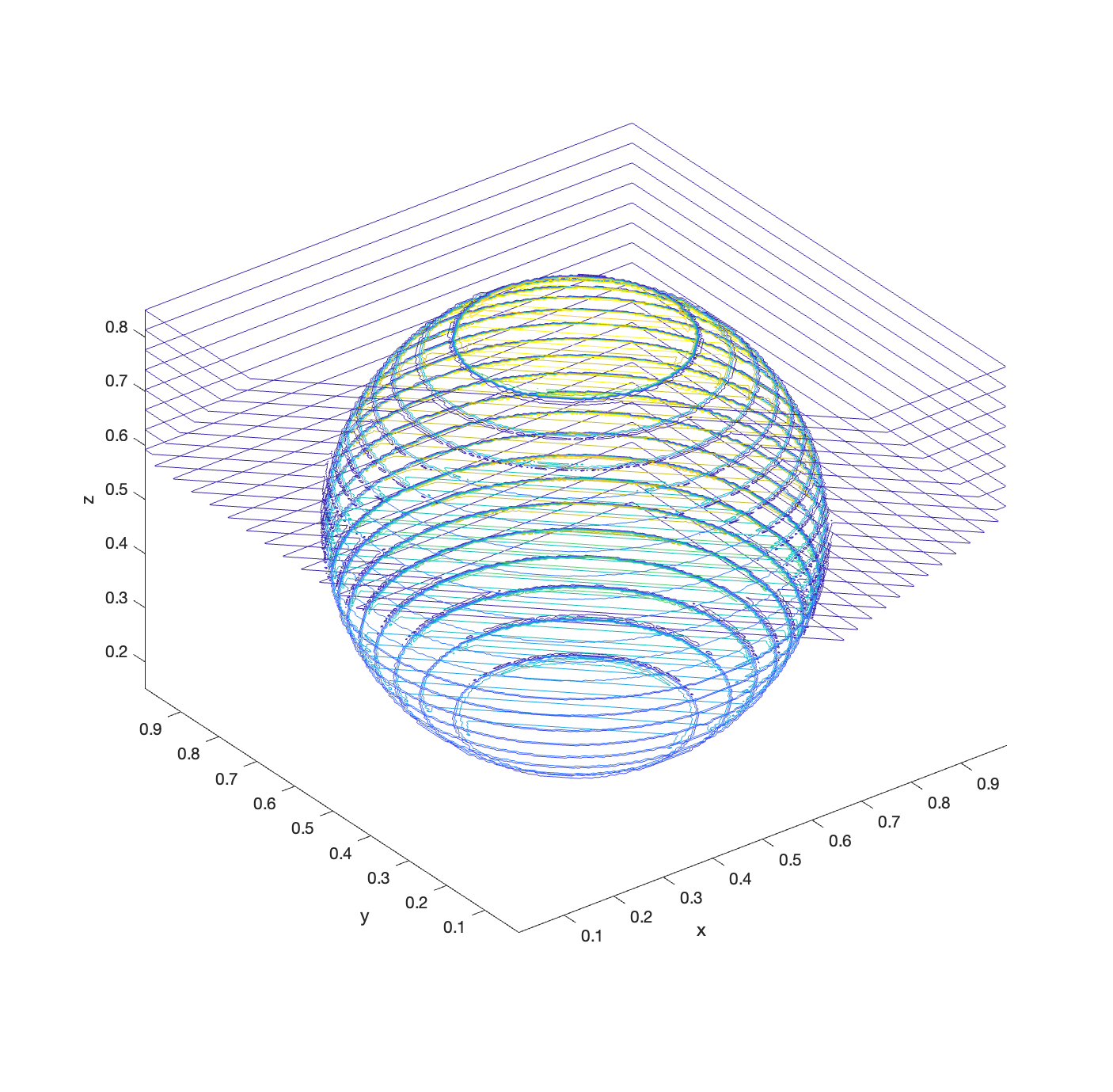,height=8cm}
\psfig{figure=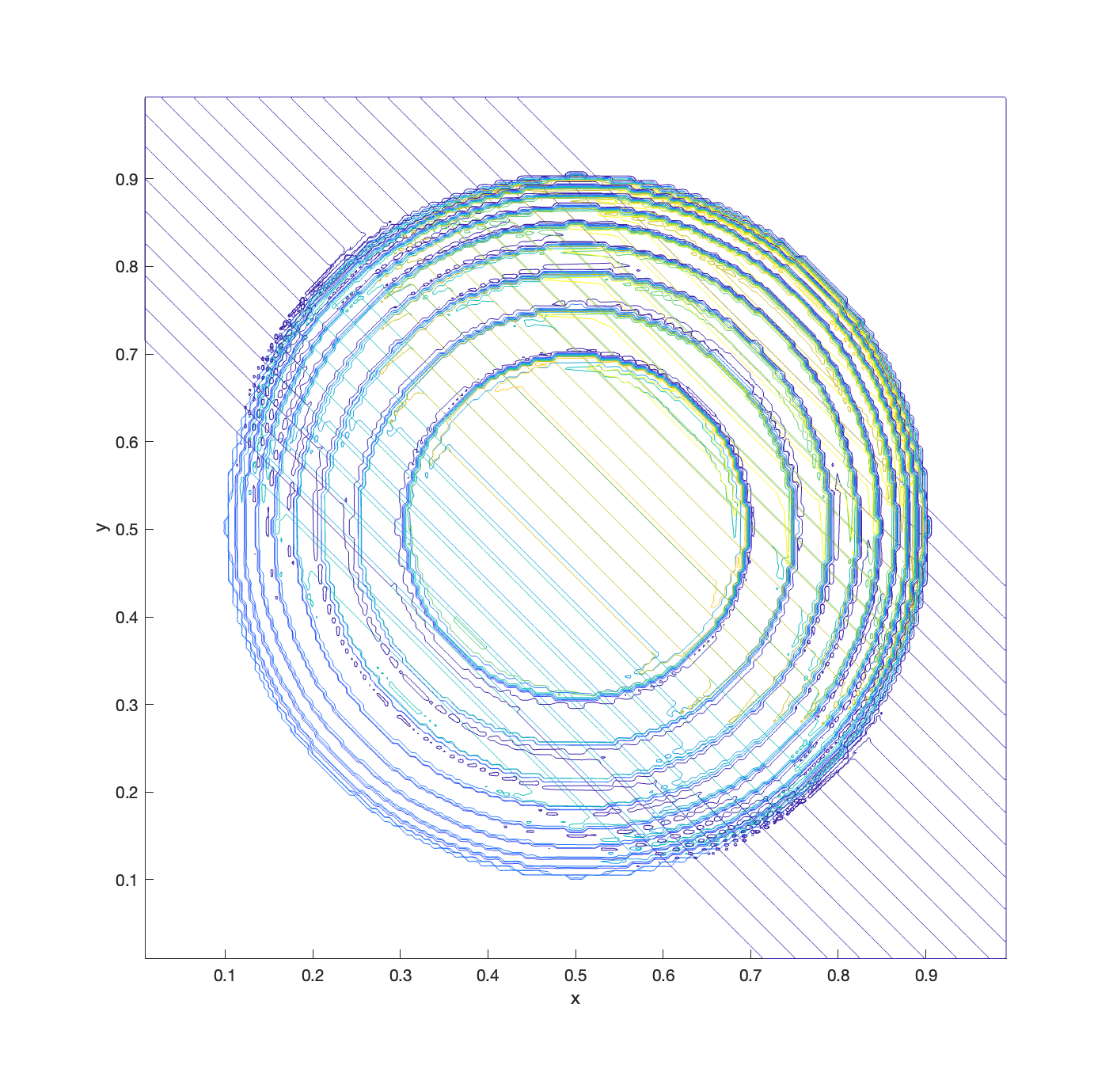,height=8cm}}
\caption{Left column, approximation of volumetric data obtained from the trivariate function in (\ref{3D}). The approximation shown has been obtained for $\psi_d$ (first row), and the linear algorithm (second row). Right column, the same figure but projected over the plane $xy$.}\label{figuras_3D}
\end{figure}

\section{Conclusions}

In this article, we have introduced a new WENO B-spline-based algorithm, where the adaption properties are reached through the nonlinear modification of the b-spline functions, which are known to be a partition of unity. These B-splines are used to construct nonlinear weight functions, that must satisfy certain properties in order to assure that the spline can avoid Gibbs-like oscillations close to discontinuities in the function. This approach is completely new and, up to our knowledge, has not been proposed in the literature before. Apart from the adaption to the discontinuities, the new algorithm conserves the smoothness of the classical spline. We have provided theoretical proofs for the accuracy of this new construction at smooth zones and close to the discontinuities, as well as conditions for the weight functions to conserve the accuracy of the data at smooth zones. Through a tensor product strategy, we are capable to extend our results to data of any number of dimensions. All the numerical experiments that we have presented  satisfy the theoretical conclusions that we have reached.

%\end{thebibliography}

\end{document}